%
%
%
%
\magnification=1200
\input amstex
\input amsppt.sty
\def\SBIMSMark#1#2#3{
 \font\SBF=cmss10 at 10 true pt
 \font\SBI=cmssi10 at 10 true pt
 \setbox0=\hbox{\SBF Stony Brook IMS Preprint \##1}
 \setbox2=\hbox to \wd0{\hfil \SBI #2}
 \setbox4=\hbox to \wd0{\hfil \SBI #3}
 \setbox6=\hbox to \wd0{\hss
             \vbox{\hsize=\wd0 \parskip=0pt \baselineskip=10 true pt
                   \copy0 \break%
                   \copy2 \break%
                   \copy4 \break}}
 \dimen0=\ht6   \advance\dimen0 by \vsize \advance\dimen0 by 8 true pt
                \advance\dimen0 by -\pagetotal
 \dimen2=\hsize \advance\dimen2 by .25 true in
%
%
  \openin2=publishd.tex
  \ifeof2\setbox0=\hbox to 0pt{}
  \else 
     \setbox0=\hbox to 3.1 true in{
                \vbox to \ht6{\hsize=3 true in \parskip=0pt  \noindent  
                \input publishd.tex 
                \vfill}}
  \fi
  \closein2
  \ht0=0pt \dp0=0pt
 \ht6=0pt \dp6=0pt
 \setbox8=\vbox to \dimen0{\vfill \hbox to \dimen2{\copy0 \hss \copy6}}
 \ht8=0pt \dp8=0pt \wd8=0pt
 \copy8
 \message{*** Stony Brook IMS Preprint #1, #2 ***}
}

\UseAMSsymbols
\loadbold
\input psbox
\psfordvips
\let\fillinggrid=\relax
\documentstyle{amsppt}
\refstyle{C}
\NoBlackBoxes
\NoRunningHeads
\vsize=21.2 truecm
\hsize=15 truecm
\hcorrection{0.8truecm}
\def\today{November, 1997}
\newcount\refno\refno=1
\define\referencelabel{\number\refno} 
\define\newrefno{\advance\refno by1}
\edef\refnoBi
{\referencelabel}\newrefno
\define\refBi{{\bf[\refnoBi]}}
\edef\refnoCa
{\referencelabel}\newrefno
\define\refCa{{\bf[\refnoCa]}}
\edef\refnodF
{\referencelabel}\newrefno

\edef\refnodFdM
{\referencelabel}\newrefno
\define\refdFdM{{\bf[\refnodFdM]}}
\edef\refnoFr
{\referencelabel}\newrefno
\define\refFr{{\bf[\refnoFr]}}
\edef\refnoGd
{\referencelabel}\newrefno
\define\refGd{{\bf[\refnoGd]}}
\edef\refnoGk
{\referencelabel}\newrefno

\edef\refnoGkone
{\referencelabel}\newrefno

\edef\refnoGS
{\referencelabel}\newrefno
\define\refGS{{\bf[\refnoGS]}}
\edef\refnoHea
{\referencelabel}\newrefno
\define\refHea{{\bf[\refnoHea]}}
\edef\refnoHeb
{\referencelabel}\newrefno
\define\refHeb{{\bf[\refnoHeb]}}
\edef\refnoLa
{\referencelabel}\newrefno
\define\refLa{{\bf[\refnoLa]}}
\edef\refnoLy
{\referencelabel}\newrefno

\edef\refnoMc
{\referencelabel}\newrefno

\edef\refnoMS
{\referencelabel}\newrefno
\define\refMS{{\bf[\refnoMS]}}
\edef\refnoRa
{\referencelabel}\newrefno
\define\refRa{{\bf[\refnoRa]}}
\edef\refnoSwa
{\referencelabel}\newrefno
\define\refSwa{{\bf[\refnoSwa]}}
\edef\refnoSwb
{\referencelabel}\newrefno
\define\refSwb{{\bf[\refnoSwb]}}
\edef\refnoYoa
{\referencelabel}\newrefno
\define\refYoa{{\bf[\refnoYoa]}}
\edef\refnoYob
{\referencelabel}\newrefno
\define\refYob{{\bf[\refnoYob]}}

%
%

\define\AnM{Ann. of Math.}

\define\CMP{Commun. Math. Phys.}
\define\ETDS{Ergod. Th. \& Dynam. Sys.}
\define\IHES{Publ. Math. IHES}

\define\JAM{J. Anal. Math.}

\define\PAMS{Proc. Amer. Math. Soc.}
\define\PCPS{Proc. Cambridge Phil. Soc.}

\define\NL{Nonlinearity}

%
%

\define\Bi{P.~Billingsley}
\define\Ca{L.~Carleson}
\define\dF{E.~de~Faria}

\define\Fr{J.~Franks}

\define\Gd{I.~J.~Good}
\define\Gk{J.~Graczyk}
\define\He{M.~Herman}
\define\La{O.~Lanford}
\define\Ly{M.~Lyubich}
\define\Mc{C.~McMullen}
\define\Me{W.~de~Melo}
\define\Ra{D.~Rand}

\define\Sw{G.~\'Swi\c{a}tek}

\define\vS{S.~van~Strien}
\define\Yo{J.~C.~Yoccoz}
\overfullrule=0pt
\def\d{\displaystyle}

\def\bk(#1){{\boldkey #1}}
\define\({\left(}
\define\){\right)}  
\define\[{\left[}
\define\]{\right]}

\let\epsilon=\varepsilon
\def\Myboxa(#1)(#2)(#3){{\hsize #1 \eightpoint{\noindent\strut
    \centerline{#2}\newline\centerline{#3}\strut}}}     
\def\Myboxb(#1)(#2)(#3)(#4){{\hsize #1 \eightpoint{\noindent\strut
    \centerline{#2}\newline\centerline{#3}\newline\centerline{#4}\strut}}}
\def\Myboxc(#1)(#2)(#3)(#4)(#5)(#6){{\hsize #1
    \eightpoint{\noindent\strut\centerline{#2}
      \newline\centerline{#3}\newline\centerline{#4}
      \newline\centerline{#5}\newline\centerline{#6} \strut}}}       
\newcount\eqnocounter
\eqnocounter=0
\define\slabel#1{\global\advance\eqnocounter by1%
\global\expandafter\edef\csname eqno#1\endcsname{(\number\eqnocounter)}%
\eqno({\number\eqnocounter})}
\define\alabel#1{\global\advance\eqnocounter by1%
\global\expandafter\edef\csname eqno#1\endcsname{(\number\eqnocounter)}%
({\number\eqnocounter})}
\define\endclaim{\egroup\par\ifdim\lastskip<\smallskipamount
\removelastskip\penalty20\smallskip\fi}
\define\claim{\smallbreak\noindent  
{\it Claim.\enspace}\ignorespaces\bgroup\sl}
\define\endassertion{\egroup\par\ifdim\lastskip<\smallskipamount
\removelastskip\penalty20\smallskip\fi}
\define\assertion(#1){\smallbreak\noindent 
{\rm(}{#1}{\rm)}\enspace\ignorespaces\bgroup\sl}
\define\endcase{\par\ifdim\lastskip<\smallskipamount
\removelastskip\penalty3000\smallskip\fi}
\define\case#1.{\smallbreak\noindent 
{\sl Case\/} #1.\enspace\ignorespaces}
\define\endcasebody{\par\ifdim\lastskip<\smallskipamount
\removelastskip\penalty-20\smallskip\fi}

%

\newcount\statemno\statemno=1
\define\newstatemno{\advance\statemno by1}
\nologo
\SBIMSMark{1997/16}{November 1997}{}
\topmatter

\title
Rigidity of critical circle mappings I
\endtitle

\author
Edson~de~Faria and Welington~de~Melo
\endauthor

\rightheadtext{Renormalization: compactness and convergence}

\abstract
We prove that two $C^r$ critical circle maps with the same rotation number of
{\it bounded type} are $C^{1+\alpha}$ conjugate for some $\alpha>0$ provided
their successive renormalizations converge together at an exponential
rate in the $C^0$ sense. The number $\alpha$ depends only on the rate of
convergence. We also give examples of $C^\infty$ critical circle maps with the
same rotation number that are not $C^{1+\beta}$ conjugate for any $\beta>0$. 
\endabstract

\date
\today
\enddate

\address
Instituto de Matem\'atica e Estat\'\i stica,
Universidade de S\~ao Paulo, \endgraf
Rua do Mat\~ao 1010, Butant\~a 
CEP05508-900 S\~ao Paulo SP - Brasil\endaddress
\email edson\@ime.usp.br\endemail

\address
Instituto de Matem\'atica Pura e Aplicada,
Estrada Dona Castorina 110,\endgraf
Jardim Bot\^anico, CEP22460
Rio de Janeiro RJ - Brasil
\endaddress
\email demelo\@impa.br\endemail


\thanks This work has been partially supported by the Pronex Project on
Dynamical Systems, \endgraf by FAPESP Grant 95/3187-4 and by CNPq Grant
30.1244/86-3.   
\endthanks


\endtopmatter

\heading  1. 
Introduction
\endheading
The purpose of this paper is to study certain rigidity questions concerning
cri\-ti\-cal circle mappings. This study is continued in {\refdFdM}.

In the qualitative theory of smooth dynamical systems, the notions of {\it
rigidity} and {\it flexibility} play an important role.
The smooth systems are usually classified according to the equivalence 
relation given by topological conjugacies: two smooth maps $f$ and $g$ are
topologically equivalent if there exists a homeomorphism $h$ of the ambient
space such that $h \circ f=g\circ h$.  
Such a homeomorphism maps orbits of $f$ onto orbits of $g$. One can also 
consider a stronger equivalence relation given by smooth conjugacies.
This leads to a quantitative or geometric classification of smooth
dynamical systems, since a smooth conjugacy, being essentially affine at small
scales, preserves the small-scale geometric properties of the dynamics.
Hence each topological equivalence class is ``foliated'' by the smooth
conjugacy classes and the quotient space is the {\it moduli} or {\it
deformation space} of the dynamics.
The moduli space describes the {\it flexibility} of the dynamics.  
When this space reduces to a single point, we are in the 
presence of {\it rigidity}. 

In general, since eigenvalues at the periodic points are smooth conjugacy 
invariants, we can hope to find rigidity only in the absence of
periodic points. From this viewpoint, the simplest case to consider is
that of circle diffeomorphisms. If $f$ is a circle diffeomorphism without
periodic points then $f$ is combinatorially equivalent to a rigid
rotation $R_\rho: x \mapsto x+\rho (\text {  mod  }1)$, in the sense that
for each $N$, the first $N$ elements of an orbit of $f$ are ordered in
the circle in the same way as the first $N$ elements of an orbit of
$R_\rho$. From Denjoy's theorem it follows that if $f$ is at least $C^2$
(or $C^1$ and its derivative has  bounded variation) then $f$ is topologically
conjugate to $R_\rho$. By a fundamental result of Herman {\refHea}, improved
by Yoccoz {\refYoa}, we have that if the rotation number $\rho$ satisfies a
Diophantine condition such as   
$$
\left|\rho-\frac pq\right|\geq \frac{C}{q^{2+\beta}}
\ , 
$$
for all rationals $p/q$, with $C>0$ and $0<\beta<1$, and if $f$ is
$C^r$, $r\geq 3$, then the conjugacy is $C^1$ (in fact it is
$C^{r-1-\beta-\epsilon}$ for every $\epsilon>0$). On the other hand,
Arnold proved that some such condition on the rotation number is essential:
there exist real analytic circle diffeomorphisms with irrational
rotation number such that the conjugacy with a rigid rotation is not
even absolutely continuous with respect to Lebesgue measure. 

Maps with periodic points cannot be rigid, but we can analyze the
rigidity of some relevant invariant set, such as an attractor of the map.
This is the situation studied by Sullivan and McMullen in the context
of unimodal maps of the interval. They considered the so-called
infinitely renormalizable maps of bounded combinatorial type. For such
maps, almost all orbits are asymptotic to a Cantor set which is the
closure of the critical orbit. They proved that if two such maps are
smooth enough and  have the same combinatorics then there exists a
$C^{1+\alpha}$ diffeomorphism of the real line that  conjugates the
restriction of the maps to the corresponding Cantor attractors. The
tools they developed have been of fundamental importance for the proof
of our results.  

Perhaps the most famous rigidity result in Geometry is 
the celebrated {\it Mostow rigidity theorem}. A special case of this
theorem states that two compact hyperbolic manifolds of dimension at
least 3 which have the same homotopy type are in fact isometric. Here
a hyperbolic manifold is the quotient space ${\Bbb H}^n/\Gamma$ of the
hyperbolic space ${\Bbb H}^n$  by a discrete group $\Gamma$ of
isometries. The hypothesis of the theorem implies the existence of a
quasiconformal homeomorphism of the sphere at infinity that
conjugates the actions of the two groups there. Such {\it a-priori}
step may be regarded as a {\it pre-rigidity} result. The rigidity is
then obtained by proving that this qc-homeomorphism is in fact
conformal, {\it i.e.\/} a Moebius transformation. 

The situation for critical circle mappings fits perfectly into this
framework. 
A critical circle mapping is a homeomorphism $f : S^1 \to S^1$
that is of class $C^r$, $r\geq 3$, and has a unique critical point
$c$ around which, in some $C^r$ coordinate system, $f$ has the form $x
\mapsto x^p$, where $p\geq 3$ is an odd integer called the {\it power
law} of~$f$. Yoccoz proved in {\refYob} that a critical circle mapping
without periodic points is topologically conjugate to an irrational
rotation. Later, in an unpublished work, he proved that the conjugacy
between two critical circle mappings with the same rotation number is
in fact quasisymmetric, {\it i.e.} there exists a constant $K\ge 1$
such that, for all pairs of adjacent intervals $I_1,I_2$ of equal
length $|I_1|=|I_2|$, we have 
$$
\frac 1K \le\frac {|h(I_1)|}{|h(I_2)|}\le K
\ .
$$
This is in contrast with the
diffeomorphism case where, without restriction on the rotation
number, the conjugacy may fail to be quasisymmetric (see {\refMS},
p.~75). 
Yoccoz's result, whose proof we present in \S 4 and Appendix B, is the exact
analogue of the pre-rigidity step in the proof of Mostow's theorem. 

\proclaim{Rigidity Conjecture} If $f, g$ are $C^3$ critical circle mappings
with the same irrational rotation number of bounded type and the same
power-law at the critical point, then there exists a $C^{1+\alpha}$ conjugacy
$h$ between $f$ and $g$ for some universal $\alpha >0$.   
\endproclaim

So far we have succeeded in proving this conjecture only when the maps are
real-analytic. Our proof involves real techniques developed in this paper, and
deformation of complex structures, developed in the next paper. 

\subhead 
1.1 Summary of results 
\endsubhead
We now present a quick summary of our results. As already
mentioned, we prove two main new theorems concerning critical circle
homeomorphisms. 

The first theorem brings forth the connection between
renormalization and rigidity in the context of circle maps. The proof is given
in \S 4.4.

\proclaim{First Main Theorem} Let $f$ and $g$ be topologically
conjugate $C^3$ critical circle maps, and let $h$ be the conjugacy
between $f$ and $g$ which maps the critical point of $f$ to the
critical point of $g$. If the partial quotients of their common rotation
number are bounded, and if their renormalizations converge together
exponentially fast in the $C^0$-topology, then $h$ is $C^{1+\alpha}$
for some $\alpha>0$.
\endproclaim

The second theorem shows that the bounded type hypothesis in the 
Rigidity Conjecture stated above cannot be removed. The proof occupies \S 5 in
its entirety.

\proclaim{Second Main Theorem} There exists an uncountable set $S$ of rotation
numbers such that for any $\rho\in S$ there exist $C^\infty$ critical circle
maps $f$ and $g$ with rotation number $\rho$ with the property that the
conjugacy between $f$ and $g$ sending the critical point of $f$ to the
critical point of $g$ is not $C^{1+\beta}$ for any $\beta>0$.
\endproclaim

The set $S$ is very small: its Hausdorff dimension is not greater than $1/2$.
But it does contain Diophantine numbers, in somewhat remarkable contrast with
the case of circle diffeomorphisms. The saddle-node surgery procedure
we develop here is quite general, and can be used to produce
similar counterexamples to the rigidity of infinitely renormalizable
unimodal maps with special unbounded combinatorics.

All estimates performed in this paper rely heavily on the {\it real a-priori
bounds} of M.~Herman {\refHeb} and G.~\'Swi\c{a}tek {\refSwa}. These bounds are
revisited in \S 3. Several technical consequences of the real bounds needed in
this paper, such as the $C^{r-1}$ boundedness of the renormalizations of a
$C^r$ critical circle map, are gathered in Appendix A.  

\subhead Acknowledgments \endsubhead
We are grateful to Dennis Sullivan for suggesting many interesting ideas, and
to Curt McMullen for explaining to us various aspects of his beautiful work on
renormalization and rigidity. We also want to thank Jean-Christophe Yoccoz for
providing us with a sketch of the proof of his unpublished theorem on
quasisymmetric conjugacies between critical circle maps. We are indebted also
to Marco Martens, Jacek Graczyk, Oscar Lanford, Rafael de la Llave, Michael
Lyubich and Henri Epstein, for several useful discussions on these and related
matters. This work has been financially sponsored by PRONEX, CNPq
and FAPESP. We also would like to thank IHES, MSRI, ETH-Z\"urich, IMPA and
SUNY at Stony Brook for their warm hospitality and generous support.

\heading 2. Preliminaries
\endheading

We have three goals in this section. First, we present some of the
basic notations commonly used when studying circle maps. Second, we
present the notions of commuting pair and renormalization in the context of
circle maps, and discuss their relationship. Third, we state the distortion
tools that are necessary for proving the real bounds in \S 3.

\subhead 
2.1 Critical circle mappings
\endsubhead
Following the tradition in this subject, we identify the unit circle
$S^1$ with the one-dimensional torus ${\Bbb R}/{\Bbb Z}$. The obvious
advantage of such identification is that it allows us to use additive
notation when dealing with circle mappings.

We briefly recall some standard facts concerning circle mappings.
Given a homeomorphism $f: S^1\to S^1$, we denote its rotation number
by $\rho(f)$. It can be expressed as a continued fraction
$$
\rho(f)\;=\;[a_0,a_1,\ldots,a_n,\ldots ]\;=\;
\cfrac 1\\ 
a_0+\cfrac 1\\ 
a_1+\cfrac 1\\ 
{}\cfrac {\cdots}\\
a_n+\cfrac 1\\ 
{\cdots}
\endcfrac
\ ,
$$
which can be finite or infinite, depending on whether $\rho(f)$ is
rational or irrational, respectively. The positive integers $a_n$ are
the {\it partial quotients} of $\rho(f)$. They give rise to a sequence of
{\it return times} for $f$, recursively defined by $q_0=1$, $q_1=a_0$
and $q_{n+1}=a_nq_n+q_{n-1}$ for all $n\ge 1$ (for which $a_n$
exists -- an assumption that will be implicit henceforth).  
Given $x\in S^1$ and $n\ge 1$, we denote by $J_n(x)$ the closed interval
containing $x$ whose endpoints are $f^{q_n}(x)$ and $f^{q_{n-1}}(x)$.
We also let $I_{n-1}(x)\subseteq J_n(x)$ be the closed interval whose
endpoints are $x$ and $f^{q_{n-1}}(x)$. Observe that $J_n(x)=I_n(x)\cup
I_{n-1}(x)$ for all $n\ge 1$. 

From the dynamics standpoint, we are not interested in {\it all}
circle homeomorphisms, but only in those that possess a unique critical
point in $S^1$, being local diffeomorphisms everywhere else. These
are the so-called {\it critical circle maps}. More precisely, let
$f:S^1\to S^1$ be a $C^r$ homeomorphism, for some $r\ge 1$. We say
that $f$ is a critical circle map if there exists $c\in S^1$ (the
critical point) such that $f'(c)=0$ and $f'(x)\neq 0$ for all $x\neq
c$. Moreover, we require $f$ to have a {\it power-law} at $c$. This
means that in a suitable $C^r$ coordinate system around the critical
point, our $f$ is represented by a map of the form $x\mapsto
x|x|^{p-1}+a$, for some real number $p>1$ called the {\it power-law
exponent} of $f$. The class of all $C^r$ critical circle maps will be
denoted by ${\roman{Crit}}^r(S^1)$. 

Since the critical point $c$ of a critical circle map is a
distinguished point on the circle, we will write $I_n$ and $J_n$
throughout, instead of $I_n(c)$ and $J_n(c)$, respectively.

\subhead
2.2 Commuting pairs 
\endsubhead
We will study the successive renormalizations of a
critical circle map $f$. Here, as in many other settings in dynamics,
the word renormalization is taken to mean a (suitably normalized)
Poincar\'e first return map of $f$ to some neighborhood of its
critical point. Abstracting the essential features of such first
return maps yields the notion of {\it commuting pair}, due to
O.~Lanford {\refLa} and D.~Rand {\refRa}. We formulate this notion as
follows.  

\demo{Definition} A {\it $C^r$ commuting pair} consists of two
mappings $\bk(f)_{-}:[\lambda,0]\to {\Bbb R}$, where $\lambda<0$,
and $\bk(f)_{+}:[0,1]\to {\Bbb R}$, satisfying the following
conditions. 
\itemitem{[P$_1$]} Both $\bk(f)_{-}$ and $\bk(f)_{+}$ are $C^r$
orientation-preserving homeomorphisms onto their images. 
\itemitem{[P$_2$]} We have $\bk(f)_{-}(0)=1$, ${\boldkey
f}_{+}(0)=\lambda$ and $0<\bk(f)_{-}(\lambda)=\bk(f)_{+}(1)<1$.
\itemitem{[P$_3$]} We have $D\bk(f)_{-}(x)>0$ for all $\lambda\le
x<0$, and $D\bk(f)_{+}(x)>0$ for all $0<x\le 1$.
\itemitem{[P$_4$]} For each $1\le k\le r$, we have 
$D^k(\bk(f)_{+}\circ \bk(f)_{-})(0)= D^k(\bk(f)_{-}\circ
\bk(f)_{+})(0)$. 

\noindent A {\it critical} commuting pair is a commuting pair such that
$D\bk(f)_{-}(0)= 0=D\bk(f)_{+}(0)$.
\enddemo

Although it is more customary to use the symbols ${\bold \xi}$ and ${\bold
\eta}$  instead of $\bk(f)_{-}$ and $\bk(f)_{+}$,
respectively, the present notation will be more convenient for our
purposes in this paper.  It can be proved that, in the presence of the
other conditions, P$_4$ is equivalent to the following. 
\itemitem{[P$_4'$]} There exist open intervals $\Delta_{-}\supseteq
[\lambda,0]$ and $\Delta_{+}\supseteq [0,1]$, and $C^r$ homeomorphic
extensions ${\bk(F)}_{-}:\Delta_{-}\to {\Bbb R}$ and
${\bk(F)}_{+}:\Delta_{+}\to {\Bbb R}$ of $\bk(f)_{-}$ and
$\bk(f)_{+}$ respectively, satisfying ${\bk(F)}_{+}\circ
{\bk(F)}_{-}(x)= {\bk(F)}_{-}\circ
{\bk(F)}_{+}(x)$ for all $x\in \Delta_{-}\cap\Delta_{+}$ 
such that ${\bk(F)}_{\pm}(x)\in \Delta_{\mp}$ (the set of such
$x$ is an open interval around $0$).

\noindent This justifies the name {\it commuting pair}. The class of
all $C^r$ critical commuting pairs will be denoted by ${\bold P}^r$.
We shall henceforth identify a commuting pair
$(\bk(f)_{-},\bk(f)_{+})$ with a single map $\bk(f):[\lambda,1]\to
[\lambda, 1]$, called the {\it shadow} of the commuting pair, defined
as follows, 
$$
\bk(f)(x)\;=\;
\cases
\bk(f)_{-}(x),&\text{when $\lambda\le x\le 0$}\cr
{}&{}\cr
\bk(f)_{+}(x),&\text{when $0< x\le 1$}.\cr
\endcases
$$
To each commuting pair $\bk(f)$ we associate an element $a\in {\Bbb
N}\cup\{\infty\}$ called the {\it height} of $\bk(f)$, in the following
way. If there exists $n\ge 1$ such that $\bk(f)^{n+1}(1)<0\le
\bk(f)^n(1)$, then we set $a=n$; otherwise we set $a=\infty$. It is
clear that $\bk(f)$ has infinite height if and only if there exists
$0<\overline{x}<1$ such that $\bk(f)(\overline{x})=\overline{x}$.

\subhead
2.3 Renormalizing a commuting pair
\endsubhead
Every commuting pair $\bk(f)$ with finite height $a$ such that
$\bk(f)^a(1)>0$ can be {\it renormalized}, in the following sense. Let
$\Lambda:{\Bbb R}\to{\Bbb R}$ be the linear map $x\mapsto \lambda x$,
let $\lambda'=\bk(f)^a(1)/\lambda<0$, and let ${\Cal
R}\bk(f):[\lambda',1]\to [\lambda',1]$ be the map defined by
$$
{\Cal R}\bk(f)(x)\;=\;
\cases
\Lambda^{-1}\circ\bk(f)\circ\Lambda(x),&\text{when $\lambda'\le x\le 0$}\cr
{}&{}\cr
\Lambda^{-1}\circ\bk(f)^{a+1}\circ\Lambda(x),&\text{when $0< x\le 1$}.\cr
\endcases
$$
This map is (the shadow of) a commuting pair $({\Cal
R}\bk(f)_{-},{\Cal R}\bk(f)_{+})$, called the {\it first
renormalization} of $\bk(f)$. Equivalently,
$$
\cases
{\Cal R}\bk(f)_{-}(x)\;=\;
\Lambda^{-1}\circ\bk(f)_{+}\circ\Lambda(x),&\text{for all $\lambda'\le
x\le 0$}\cr 
{}&{}\cr
{\Cal R}\bk(f)_{+}(x)\;=\;
\Lambda^{-1}\circ\bk(f)_{+}^{a}\circ\bk(f)_{-}\circ\Lambda(x),&\text{for
all $0\le x\le 1$}.\cr 
\endcases
$$
The class of all $C^r$ critical commuting pairs which are
renormalizable in this sense will be denoted ${\bold P}^r_1$. In this
way, we have a well-defined map ${\Cal R}:{\bold P}^r_1\to {\bold
P}^r$, the so-called {\it renormalization operator}. More generally,
for all $n\ge 1$ we write ${\bold P}^r_n={\Cal R}^{-n}({\bold P}^r)$
for the set of all $C^r$ critical commuting pairs which can be
renormalized $n$ times. We have ${\bold P}^r_{n+1}\subseteq {\bold
P}^r_n$ for all $n$. We are especially interested in the set of all
{\it infinitely renormalizable} critical commuting pairs, namely
$$
{\bold P}^r_{\infty}\;=\;\bigcap_{n\ge 1} {\bold P}^r_n
\ .
$$
Given $\bk(f)\in {\bold P}^r$, let $a_0=a$ be its height, and for each
$n\ge 1$ such that $\bk(f)\in {\bold P}^r_n$, let $a_n$ be the height
of ${\Cal R}^n(\bk(f))$. This can be a finite or infinite sequence; in
any case, using the convention $1/\infty=0$, we define the {\it
rotation number} of $\bk(f)$ to be 
$$
\rho(\bk(f))\;=\;[a_0,a_1,\ldots,a_n,\ldots ]\;=\;
\cfrac 1\\ 
a_0+\cfrac 1\\ 
a_1+\cfrac 1\\ 
{}\cfrac {\cdots}\\
a_n+\cfrac 1\\ 
{\cdots}
\endcfrac
\ .
$$
In particular, $\rho({\Cal R}\bk(f))=[a_1,a_2,\ldots ]$, that is, the
renormalization operator acts as the Gaussian shift on continued
fractions. 
\bigskip

$$
\psannotate{
\psboxto(8.0cm;0cm){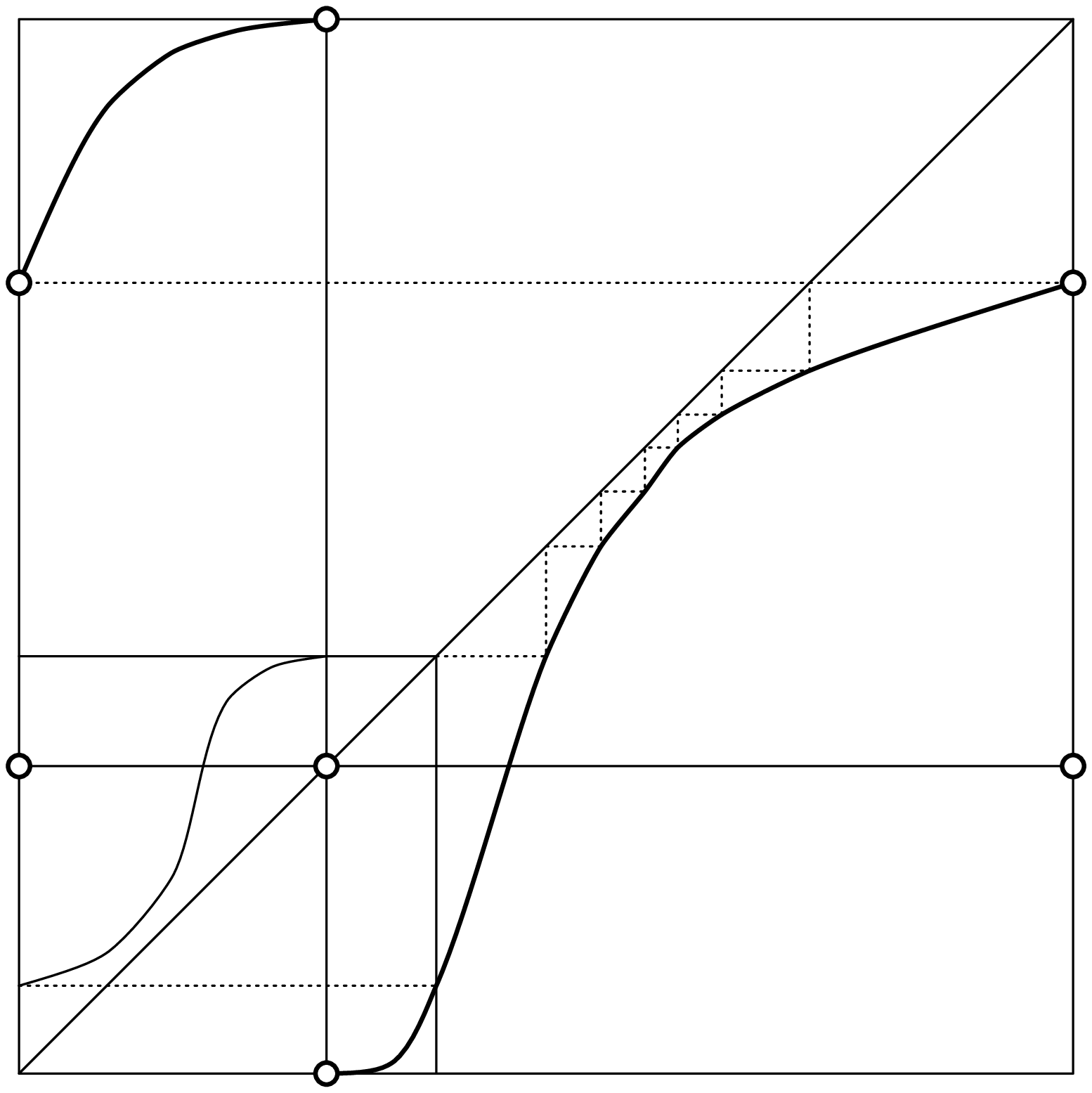}}
{\fillinggrid\at{-1.05\pscm}{4.05\pscm}{$f^{q_n}(c)$}
             \at{5.4\pscm}{4.05\pscm}{$c$}
             \at{16.1\pscm}{4.05\pscm}{$f^{q_{n-1}}(c)$}
             \at{2.8\pscm}{14\pscm}{$f^{q_{n-1}}$}
             \at{13.5\pscm}{10\pscm}{$f^{q_n}$}}
$$
\bigskip
\centerline{\sl Figure 1. Two consecutive renormalizations of $f$.}
\bigskip

\subhead
2.4 Renormalizing a critical circle map
\endsubhead
Let $f$ be a critical circle map with critical point $c$, and
for each $k\ge 0$ such that $f^{q_k}(c)\neq c$, let $A_k:{\Bbb
R}\to {\Bbb R}/{\Bbb Z}$ be the affine covering map such that
$A_k([0,1])=I_k$, with $A_k(0)=c$ and $A_k(1)=f^{q_k}(c)$. For
each $n\ge 1$ such that $f^{q_k}(c)\neq c$ for all $0\le k\le n$,
consider the Poincar\'e first return map $f_n:I_n\cup I_{n-1}\to
I_n\cup I_{n-1}$, namely 
$$
f_n(x)\;=\;\cases
f^q(x)&\text{when $x\in I_n$}\cr
{}&{}\cr
f^Q(x)&\text{when $x\in I_{n-1}$},\cr
\endcases
$$
where $q=q_{n-1}$ and $Q=q_n$. Define $\lambda_n$ to be the largest
negative number such that $A_{n-1}(\lambda_n)=f^Q(c)$ (one sees in
fact that $\lambda_n=-|I_n|/|I_{n-1}|$). 
Then $A_{n-1}([\lambda_n,0])= I_n$ and $A_{n-1}([0,1])= I_{n-1}$, and
we can consider the map $\bk(f)_n:[\lambda_n,1]\to [\lambda_n,1]$
given by $\bk(f)_n\;=\;A_{n-1}^{-1}\circ f_n\circ A_{n-1}$. Here, it
is implicit that $A_{n-1}^{-1}$ is the inverse branch that maps
$I_n\cup I_{n-1}$ onto $[\lambda_n,1]$. This defines (the shadow of) a 
$C^r$ critical commuting pair called {\it $n$-th
renormalization of $f$}. It is well-defined provided the rotation
number of $f$ has a  continued-fraction development of length at least
$n+1$ (in particular it is well-defined for all $n$ when the rotation
number of $f$ is irrational). It is easy to see in this case that
$\bk(f)_{k+1}={\Cal R}\bk(f)_k$ for all $1\le k\le n-1$. Moreover, if
$a_0,a_1,\ldots$ are the partial quotients of the rotation number of $f$
then, from the recurrence relations satisfied by the sequence of return
times $q_k$, we see at once that the height of $\bk(f)_k$ is equal to
$a_k$, and that $\rho(\bk(f)_k)=[a_k,a_{k+1},\ldots ]$.

\demo{Remark} Note that the largest interval containing $I_n$ on which
$f^{q_{n-1}}$ is a diffeomorphism is the $[\alpha_n,c]$ where
$f^{a_{n-1}q_{n-1}}(\alpha_n)=f^{q_n}(c)$. Similarly, the largest interval
containing $I_{n-1}$ on which $f^{q_n}$ is a diffeomorphism is $[c,\beta_n]$
where $f^{q_n}(\beta_n)=f^q_{n-1}(c)$.
\enddemo

\subhead
2.5 The $C^k$ metrics
\endsubhead
The following is only one of several equivalent ways of defining a $C^k$
distance between commuting pairs. We normalize our commuting pairs to be
defined on $[0,1]$, using for each $\bk(f)$ a fractional linear transformation
that maps $\lambda,0, 1$  respectively to $0,\frac{1}{2}, 1$, and then use the
$C^k$ norm of the difference of the normalized pairs. The $C^k$ distance  
between $\bk(f)$ and $\bk(g)$ defined in this fashion is 
denoted by $d_k(\bk(f),\bk(g))$.

Let us be a bit more precise. If a function $\varphi:[0,1]\to {\Bbb R}$ 
has a jump discontinuity at $x=1/2$ but is elsewhere $k$ times 
continuously differentiable, let $\|\varphi\|_k=\max\{\|\varphi^-\|_k,
\|\varphi^+\|_k\}$ where $\varphi^-$ is the restriction of $\varphi$
to $[0,\frac{1}{2}]$ and  $\varphi^+$ is the restriction of $\varphi$
to $[\frac{1}{2},1]$ . Given two elements $\bk(f):[\lambda,1]\to {\Bbb R}$ 
and $\bk(g):[\mu,1]\to {\Bbb R}$ of ${\bold P}^r$, and given $0\le
k\le r$, let $A_\lambda$ be the fractional linear transformation which
maps $0,\lambda,1$ to $0,\frac{1}{2},1$, respectively, and let $A_\mu$
be similarly defined. Then write
$$
d_k(\bk(f),\bk(g))\;=\;\max\left\{|\lambda-\mu|,\|A_\lambda \bk(f)
A_\lambda^{-1} - A_\mu \bk(g) A_\mu^{-1}\|_k\right\}
\ .
$$

\proclaim{Lemma {2.1}} For each $0\le k\le r$, $d_k$ is a metric.
\endproclaim
\demo{Proof} The only thing not entirely obvious is that
$d_k(\bk(f),\bk(g))=0$ implies $\bk(f)=\bk(g)$. But if
$d_k(\bk(f),\bk(g))=0$ then on one hand $\lambda=\mu$, so that
$A_\lambda=A_\mu$, and on the other hand $A_\lambda \bk(f)
A_\lambda^{-1} - A_\mu \bk(g) A_\mu^{-1}=0$, so that 
$\bk(f)=A_\lambda^{-1}A_\mu \bk(g) A_\mu^{-1}A_\lambda =\bk(g)$.\qed
\enddemo

\proclaim{Proposition {2.2}} Let $\bk(f):[\lambda,1]\to [\lambda,1]$
and $\bk(g):[\mu,1]\to [\mu,1]$ be elements of ${\bold P}^r_\infty$,
and suppose there exists a $C^r$ diffeomorphism $h:[\lambda,1]\to
[\mu,1]$ such that $h\circ\bk(f)=\bk(g)\circ h$. Then for all
$k\le r-1$ the distances $d_k\left({\Cal
R}^n(\bk(f)),{\Cal R}^n(\bk(g))\right)$ converge to $0$ at an
exponential rate. 
\endproclaim

\demo{Proof} Let $\bk(f)_n={\Cal R}^n\bk(f)$ and $\bk(g)_n={\Cal R}^n
\bk(g)$. Then $\bk(f)_n=h_n^{-1}\circ\bk(g)_n\circ h_n$, where $h_n$
is obtained from $h$ by restriction and affine rescaling. We will see
below (after we prove the real bounds for critical circle maps, {\it
cf.\/} Theorem {3.1}) that $\{h_n\}$ converges exponentially in the
$C^r$ sense to the space of {\it affine} maps. Therefore, we have that
$d_{r-1}(\bk(f)_n,\bk(g)_n)\to 0$ exponentially fast. \qed
\enddemo

\subhead
2.6 Distortion tools
\endsubhead
In \S 3 we will need some distortion tools to get real bounds for critical
circle maps. The most basic is the notion of {\it cross-ratio
distortion}. Given intervals $M\subseteq T$ on the line or circle, their
cross-ratio is defined as
$$
D(M,T)\;=\;\frac{|M|\,|T|}{|L|\,|R|}
\ ,
$$
where $L$ and $R$ are the left and right components of $T\setminus M$. The
cross-ratio distortion of a map $f$ (whose domain contains $T$) on the pair of
intervals $(M,T)$ is
$$
B(f;M,T)\;=\;\frac{D(f(M),f(T))}{D(M,T)}
\ .
$$
Cross-ratios are always increased by a map with negative Schwarzian
derivative. More precisely, if $f$ is $C^3$ and $Sf<0$ then $B(f;M,T)>1$.

\proclaim{Lemma {2.3}} (Cross-ratio distortion principle)

\noindent
Given a map $f$ as above, $m\ge 1$ and intervals $M\subseteq T$ such that
$f^m|T$ is a diffeomorphism onto its image, we have
$$
B(f^m;M,T)\;\ge\;\exp\{-\sigma\sum_{j=0}^{m-1}|f^j(T)|\}
\ ,
$$
where $\sigma>0$ depends on $f$ and $\max_{0\le j\le m-1}|f^j(T)|$.
\endproclaim

For a proof of (a much more general version of) this principle, see {\refMS},
p.~287. This fact will be used in combination with the following classical
distortion principle. For intervals $M\subseteq T$ as above we define the {\it
space} of $M$ inside $T$ to be the smallest of the ratios $|L|/|M|$ and
$|R|/|M|$.

\proclaim{Lemma {2.4}} (Koebe distortion principle)

\noindent
Given $\ell, \tau>0$ and a map $f$ as above, there
exists $K=K(\ell,\tau,f)>1$ of the form
$$
K\;=\;\left(1+\frac{1}{\tau}\right)^2\,\exp{C\ell}
\ ,
$$
where $C$ is a constant depending only on $f$, with the following property. If
$T$ is an interval such that $f^m|T$ is a diffeomorphism and if
$\sum_{j=0}^{m-1}|f^j(T)|\le \ell$, then for each interval $M\subseteq T$ for
which the space of $f^m(M)$ inside $f^m(T)$ is at least $\tau$ and for all 
$x,y\in M$ we have
$$
\frac{1}{K}\;\le\;\left|\frac{Df^m(x)}{Df^m(y)}\right|\;\le\;K
\ .
$$
\endproclaim

Once again, see {\refMS}, p.~295, for a proof. Used in combination with Lemma
{2.3}, the Koebe distortion principle allows one to propagate space around
under fairly general circumstances.

\heading 3.
The real a-priori bounds
\endheading

In this section we establish real {\it a-priori} bounds for critical
circle maps, obtaining as a corollary the fact that their
renormalizations are pre-compact in the $C^1$ topology. The results
are well-known, and the reader will not fail to notice the overlap
with some of the material in {\refSwb} and {\refGS}.

Let $f\colon S^1\to S^1$ be a critical circle homeomorphism with
critical point $c$.  The iterates of $c$ are denoted by 
$c_i= f^i(c)$. Let $I_n$ be the interval with endpoints $c$ and
$c_{q_n}$ that contains $c_{q_{n+2}}$, as defined in section 2.
For simplicity, we write $I_n^j=f^j(I_n)$ for all $j$ and $n$. The
most basic combinatorial fact to be remembered here is that the
collection of intervals  
$$
{\Cal P}_n\;=\;\Big\{I_{n-1},\,I_{n-1}^1,\ldots
,\, I_{n-1}^{q_n-1}\Big\}\cup\Big\{I_n,\,I_n^1,\ldots
,\,I_n^{q_{n-1}-1}\Big\} 
$$
constitutes a partition of $S^1$ modulo endpoints, called the {\it
dynamical partition of level $n$} associated to $f$. In order to get
an actual partition we exclude from each interval in ${\Cal P}_n$ its
right endpoint, say, according to the standard choice of orientation
of $S^1$.  
Let $\Cal P_n(x)$ denote the 
atom of the partition $\Cal P_n$ that contains $x$ (in particular, 
$\Cal P_n(c)$ is either $I_n$ or $I_{n-1}$ according to the parity of $n$).

\proclaim{Theorem {3.1}} (Real Bounds) 
Let $f\in {\roman{Crit}}^r(S^1)$ be a map with irrational rotation number.
There exist constants $C_0>1$ and $0<\mu_0<\mu_1<1$ depending only on
$f$ such that 
\itemitem{($a$)} If $I$ and $J$ are any two adjacent atoms of $\Cal P_n$, then
$C_0^{-1}|J|<|I|<C_0|J|$; 
\itemitem{($b$)} For every $x\in S^1$, we have $ |\Cal P_n(x)|< \mu_1 
|{\Cal P}_{n-1}(x)| $;
\itemitem{($c$)} If the rotation number of $f$ is of bounded type then
$ |\Cal P_n(x)| > \mu_0^n/C_0$;
\itemitem{($d$)} If the rotation number of $f$ is of bounded type then
$ |{\Cal P}_n(x)| > |{\Cal P}_{n-1}(x)|/C_0$;
\itemitem{($e$)} If $0<i\le j\le q_n$ then the distortion of the
restriction of $f^{j-i}$ to $I_{n-1}^i= f^i(I_{n-1})$ is bounded by $C_0$. 

\noindent In particular, the critical commuting
pairs ${\Cal R}^nf$ form a bounded sequence in the $C^1$ topology.
\endproclaim

Later in this section we will see that the bounds in this theorem are
eventually universal.

\subhead
3.1 Bounding space
\endsubhead
In what follows, two positive numbers $a$ and $b$ are said to be
{\it comparable modulo $f$}, or simply {\it comparable}, if there
exists a constant $C>1$, depending only on our map $f$, such that
$C^{-1}b\le a\le Cb$. This relation is denoted by $a\asymp b$. It is
also convenient to write $a\preccurlyeq b$ to indicate that $a\le
Cb$. Comparability modulo $f$ is reflexive and symmetric, but not
transitive since the constants multiply. Hence, if $b_1\asymp
b_2\asymp\cdots\asymp b_N$, we can only say that $b_1\asymp b_N$ if
$N$ is bounded (by a constant depending only on $f$). 
\proclaim{Lemma {3.2}} For each $n\ge 0$ there exist
$z_1,\,z_2,\,z_3,\,z_4,\,z_5\in S^1$ with $z_{j+1}=f^{q_n}(z_j)$ such
that $|z_1-z_2|\asymp |z_2-z_3| \asymp |z_3-z_4| \asymp |z_4-z_5|$.
\endproclaim
\demo{Proof} Let $z\in S^1$ be a point such that $|f^{q_n}(z)-z|\le
|f^{q_n}(x)-x|$ for all $x\in S^1$. From Koebe's principle applied
successively to $f^{-q_n}$, $f^{-2q_n}$ and $f^{-3q_n}$, we have
$$
|z-f^{q_n}(z)|\succcurlyeq |f^{-q_n}(z)-z|\succcurlyeq
|f^{-2q_n}(z)-f^{-q_n}(z)|\succcurlyeq |f^{-3q_n}(z)-f^{-2q_n}(z)|
\ .
$$
Moreover, by our choice of $z$ we have $|z-f^{q_n}(z)|
\preccurlyeq |f^{-3q_n}(z)-f^{-2q_n}(z)|$. Therefore we can take
$z_5=f^{q_n}(z),\,z_4=z,\ldots,\,z_1=f^{-3q_n}(z)$ as the desired five
points. \qed
\enddemo
\proclaim{Lemma {3.3}} Let $z_1,\,z_2,\ldots ,\,z_5$ and $w_0,\,w_1,\ldots
,\,w_5$ be points on the circle such that $z_{j+1}=f^{q_n}(z_j)$ and
$w_{j+1}=f^{q_n}(w_j)$, and such that $w_1$ lies on the interval of
endpoints $z_1$ and $z_2$ in the partition of $S^1$ determined by the
$z_i$'s. If $|z_1-z_2|\asymp |z_2-z_3| \asymp |z_3-z_4| \asymp
|z_4-z_5|$, then
$$
|w_0-w_1|\succcurlyeq |w_1-w_2|\preccurlyeq |w_2-w_3| \ .
\slabel{w}
$$
\endproclaim
\demo{Proof} Let $\ell =\min |z_j-z_{j+1}|$. If there is a $j$ with
$1\le j\le 3$ such that $|w_j-w_{j+1}|\le \ell/2$, then we must have
$|w_{j-1}-w_j|\ge \ell/2$ and $|w_{j+1}-w_{j+2}|\ge \ell/2$ also. But
then $[w_j,w_{j+1}]$ has space on both sides inside
$[w_{j-1},w_{j+2}]$. Applying $f^{-(j-1)q_n}$ to these points and
using the Koebe principle, we get {\eqnow}. If on the other hand there is
no $j$ with that property, then $|w_1-w_2|\asymp |w_2-w_3|\asymp
|w_3-w_4|$. Again, applying $f^{-q_n}$ and using Koebe we get
{\eqnow}.\qed 
\enddemo
\proclaim{Lemma {3.4}} For all $n\ge 0$ and all $x\in S^1$, we have
$|f^{q_n}(x)-x|\asymp |x-f^{-q_n}(x)|$.
\endproclaim
\demo{Proof} To show that $|x-f^{-q_n}(x)|\ge C^{-1} |f^{q_n}(x)-x|$,
let $i\le q_n$ be such that $f^i(x)\in [z_1,z_2]$, where
$z_1,\,z_2,\ldots$ are the points given by Lemma {3.2}. Then let
$w_0=f^{i-q_n}(x)$, $w_1=f^i(x)$, etc. We know from Lemma {3.3} that
$|w_0-w_1|\succcurlyeq |w_1-w_2|\preccurlyeq |w_2-w_3|$. Applying $f^{-i}$ to
these points and using the Koebe distortion principle, we find a definite
{\it space} around $[x,f^{q_n}(x)]$ inside $[f^{-q_n}(x),
f^{2q_n}(x)]$. Therefore $|x-f^{-q_n}(x)|\succcurlyeq |f^{q_n}(x)-x|$. The
proof of the opposite inequality is similar. \qed 
\enddemo
\bigskip
$$
\psannotate{
\psboxto(8.5cm;0cm){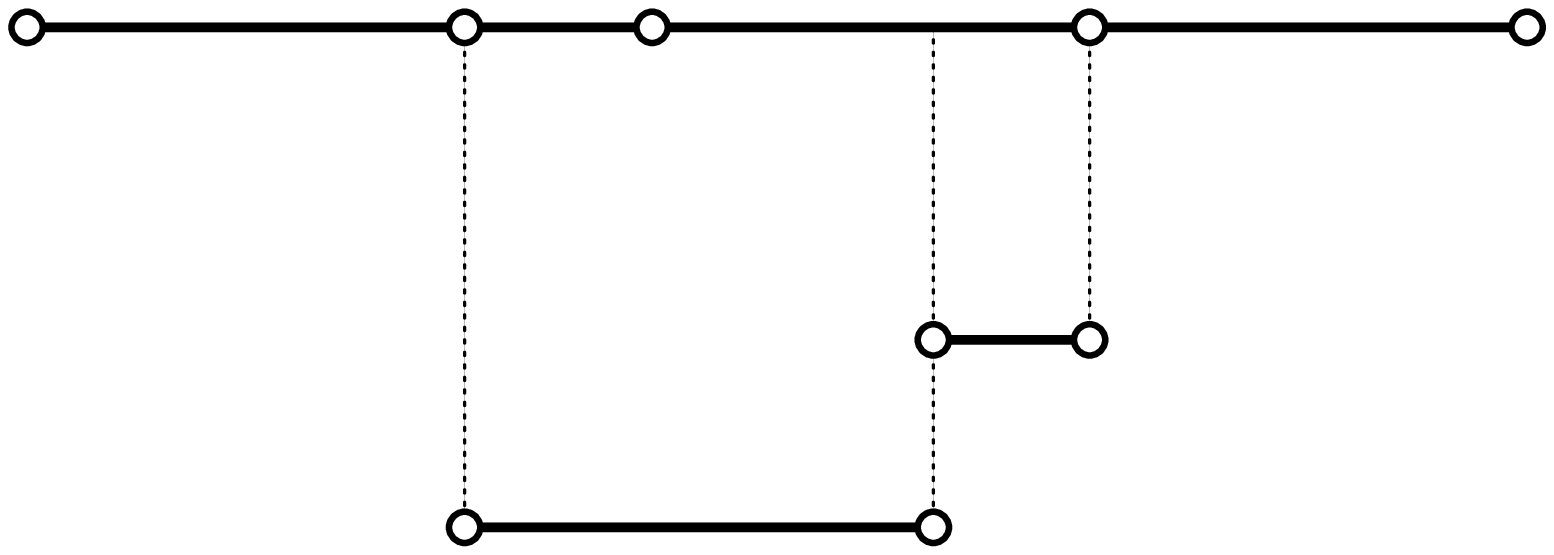}}
{\fillinggrid\at{2.6\pscm}{5.8\pscm}{$I_{n-1}^{q_n-q_{n-1}}$}
             \at{6.4\pscm}{5.8\pscm}{$I_n$}
             \at{9\pscm}{5.8\pscm}{$I_{n-1}$}
             \at{13.2\pscm}{5.8\pscm}{$I_{n-1}^{q_{n-1}}$}
             \at{7.1\pscm}{-0.7\pscm}{$I_{n-1}^{q_n}$}
             \at{10.1\pscm}{1.3\pscm}{$I_n^{q_{n-1}}$}
}
$$
\bigskip
\centerline{\sl Figure 2. These six intervals are pairwise comparable.}
\bigskip

We arrive at the following fundamental fact first proved by \'Swi\c{a}tek
{\refSwa} and Herman {\refHeb}.

\proclaim{Lemma {3.5}} Any two adjacent intervals in the dynamical
partition of level $n$ of $f$ are comparable.
\endproclaim
\demo{Proof} First we prove that all intervals in Figure 2 are pairwise 
comparable, through the following steps.

\itemitem{($a$)} From Lemma {3.4}, we know that $|I_{n-1}|\asymp
|I_{n-1}^{q_{n-1}}|$ and $|I_{n-1}^{q_n-q_{n-1}}|\asymp
|I_{n-1}^{q_n}|$. 

\itemitem{($b$)} Since the dynamical symmetric of $I_n$, namely the
interval $I_n^{-q_n}$, is contained 
in $I_{n-1}$, we also have $|I_n|\preccurlyeq |I_{n-1}|$. 

\itemitem{($c$)} Since the dynamical symmetric of $I_{n-1}$, namely
$I_{n-1}^{-q_{n-1}}$ is contained in $I_n \cup I_{n-1}^{q_n-q_{n-1}}$, we have
$|I_{n-1}|\preccurlyeq |I_{n-1}^{q_n-q_{n-1}}|$. Moreover, since
$I_{n-1}^{q_n}\subseteq I_n\cup I_{n-1}$, items ($a$) and ($b$) yield 
$|I_{n-1}^{q_n-q_{n-1}}|\preccurlyeq |I_{n-1}|$. Therefore $|I_{n-1}|\asymp
|I_{n-1}^{q_n-q_{n-1}}|$.

\itemitem{($d$)} Next, we claim that $|I_n|\asymp |I_n^{q_{n-1}}|$. To see
why, consider the diffeomorphism 
$$
f^{q_n-q_{n-1}}: I_{n-1}\cup I_{n-1}^{q_{n-1}}\to I_{n-1}^{q_n-q_{n-1}}\cup
I_{n-1}^{q_n}
\ .
$$
By the cross-ratio inequality (Lemma {2.3}) applied to $M=I_n^{q_{n-1}}$ and
$T=I_{n-1}\cup I_{n-1}^{q_{n-1}}$, we have $|I_n^{q_{n-1}}|\preccurlyeq
|I_n^{q_n}|\asymp |I_n|$. Conversely, considering the diffeomorphism
$$
f^{q_{n-1}}: I_{n-1}^{q_n-q_{n-1}}\cup I_n\to I_{n-1}^{q_n}\cup I_n^{q_{n-1}} 
$$
and applying the cross-ratio inequality to $M=I_n^{q_n}$ and
$T=I_{n-1}^{q_n-q_{n-1}}\cup I_n$, we get 
$$
|I_n|\;\asymp\;|I_n^{q_n}|\;\preccurlyeq\; |I_{n-1}^{q_n+q_{n-1}}|\;\asymp
|I_n^{q_{n-1}}|
\ .
$$
This proves our claim.

\itemitem{($e$)} Finally, we claim that $|I_{n-1}|\preccurlyeq |I_n|$, thereby
reversing the inequality in ($b$). It is here that we use the critical point
in a crucial way. Let $\theta_n=|I_n|/|I_{n-1}|$; we already know that
$\theta_n\preccurlyeq 1$. Look at the intervals $I_{n-1}^1$, $I_n^1$ and 
$I_{n-1}^{q_n-q_{n-1}+1}$, all near the critical value of $f$. By an argument
similar to the one in ($c$), we have $|I_{n-1}^1|\asymp
|I_{n-1}^{q_n-q_{n-1}+1}|$. Moreover, $|I_n^1|\asymp \theta_n^p |I_{n-1}^1|$,
where $p>1$ is the power-law of $f$ at the critical point. Hence these three
intervals have a cross-ratio comparable to $\theta_n^p$. On the other hand the
map $f^{q_{n-1}-1}$ carries them diffeomorphically onto
$I_{n-1}^{q_{n-1}}$, $I_n^{q_{n-1}}$ and $I_{n-1}^{q_n}$, respectively, whose
cross-ratio is comparable to $|I_n^{q_{n-1}}|/|I_{n-1}^{q_n}|$, which in turn
is comparable to $\theta_n$. Applying the Koebe distortion principle, we see
that $\theta_n^p \succcurlyeq\theta_n$, and so $\theta_n\succcurlyeq 1$ as
claimed. 

This proves that all six intervals in Figure 2 are comparable. To derive the
remaining comparability relations, propagate this information using Koebe's
distortion principle.\qed 
\enddemo

\demo{Proof of Theorem {3.1}} Part ($a$) is Lemma {3.5} above. The remaining
statements are straightforward consequences of ($a$). \qed
\enddemo

\subhead
3.2 Beau property of renormalization
\endsubhead
The bounds obtained in the proof of Theorem {3.1} depended on $f$, more
precisely on the space that each atom of ${\Cal P}_n$ enjoys relative to
its two neighbors in ${\Cal P}_n$. We now concentrate in proving that
such bounds eventually become universal. 
It suffices to prove that the space
in question is eventually universal. Bounds of this type are called {\it beau}
by Sullivan.  

\proclaim{Lemma {3.6}} There exists $n_0=n_0(f)$ such that for
all $n\ge n_0$ the first return map $f_n:J_n\to J_n$ satisfies
$Sf_n(x)<0$ for all $x\in J_n$. 
\endproclaim
\demo{Proof} This is proved in Theorem {A.4} of Appendix A.\qed
\enddemo

\proclaim{Lemma {3.7}} Given $\varepsilon>0$, there exists
$n_1=n_1(f,\varepsilon)>n_0(f)$ such that the following holds for all $n\ge
n_1$. 
Let $\Delta\in {\Cal P}_n$, let $k\ge 1$ be an integer such that
$f^j(\Delta)$ is contained in an element of ${\Cal P}_n$ for all $1\le j\le
k$, and let $\Delta^*$ be the union of $\Delta$ with its left and right
neighbors in ${\Cal P}_n$. 
Then we have $f^k|\Delta^*=\phi_1\circ\phi_2\circ\phi_3$ where $\phi_1$ and
$\phi_3$ are diffeomorphisms with distortion bounded by $1+\varepsilon$ and
$\phi_2$ is either the identity or a map with negative Schwarzian
derivative. In particular, if $\varepsilon$ is small enough and if
$I_{n-1}\neq \Delta\neq I_n$, then the distortion of $f^k|\Delta$ is bounded
from below by one-half. 
\endproclaim
\demo{Proof} Let $n_1>n_0$ be such that $\mu_0^{n_1-n_0}<\!\!<\varepsilon$,
where $\mu_0$ is the constant of Theorem {3.1}. 
For $n\ge n_1$, $\Delta$ and $k$ as in the statement, let
$J\in {\Cal P}_{n_0}$ be such that $\Delta\subseteq J$, let $J^*$ be the
union of $J$ with its two neighbors in ${\Cal P}_{n_0}$, and note that the
space of $\Delta^*$ inside $J^*$ is bounded from below by $C|J^*|/|\Delta^*|$,
for some constant $C>0$. Let $m\ge 0$ be the smallest integer such that
$f^m(J)\subseteq J_{n_0}$. Then for all $j\le m$ the map $f^j|J^*$ is a
diffeomorphism onto its image and, by Theorem {3.1} ($b$) and the Koebe
distortion principle, its distortion on $\Delta^*$ is bounded by 
$$
\left(1+C\frac{|\Delta^*|}{|J^*|}\right)^2
\exp\left\{C\frac{|\Delta^*|}{|J^*|}\right\}\;\le\; 
\exp\{C\mu_0^{n_1-n_0}\}\;\le\; 1+\varepsilon
\ .
$$
Now, there are two possibilities. The first is that $m\ge k$; in this
case we can take $\phi_1=f^k|\Delta^*$ and $\phi_2=\phi_3=$ identity
map. The second is that $m<k$. In this case we consider the first
return map $f_{n_0}:J_{n_0}\to J_{n_0}$ and let $\ell\ge 0$ be the
{\it largest} such that 
$$
f^k=f^{k_1}\circ f_{n_0}^\ell\circ f^{k_3}
\ ,
$$
where $k_1\ge 0$ and $k_3\ge 0$. We then take $\phi_1=f^{k_1-1}$,
$\phi_2=f\circ f_{n_0}^\ell$ and $\phi_3=f^{k_3}|\Delta^*$ (if $k_1=0$ we
take instead $\phi_2=f_{n_0}^\ell$ and $\phi_1=$ identity). By Lemma
{3.6}, $S\phi_2<0$, and by the above remarks the distortions of both
$\phi_1$ and $\phi_3$ are bounded by $1+\varepsilon$ in the
appropriate domains. \qed
\enddemo

\proclaim{Proposition {3.8}} All bounds in Theorem {3.1} are
beau. In other words, there exist universal constants $K_0>0$ and 
$0<\lambda_0<\lambda_1<1$ and some $\overline{n}=\overline{n}(f)>0$
such that for all $n\ge \overline{n}$ the constants $C_0$, $\mu_0$ and
$\mu_1$ in Theorem {3.1} can be replaced by $K_0$, $\lambda_0$ and
$\lambda_1$, respectively.
\endproclaim
\demo{Proof} This is straightforward from Lemma {3.7}.\qed
\enddemo

\demo{Remark} From now on, whenever we say that a constant
``depends only on the real bounds'', we mean that the said constant
is a universal function of constants $K_0$, $\lambda_0$ and
$\lambda_1$ given by this proposition. 
\enddemo

\head 4. How smooth is the conjugacy?
\endhead

Now we turn to the first main result in this paper. The theorem states that 
if the successive renormalizations of two critical circle maps with the same
rotation number of bounded type converge together at an exponential rate, then
such maps are $C^{1+\alpha}$ conjugate, for some universal
$\alpha>0$. First, in order to get bounds that do not depend on the maximum
of the partial quotients of the rotation number, we need to perform some
``saddle-node'' estimates and constructions.

\subhead
4.1 Saddle-node geometry 
\endsubhead
Let $a$ be a positive integer and let $\Delta_1,\Delta_2,\ldots,\Delta_{a+1}$
be consecutive intervals on the line or circle.  By an {\it almost parabolic
map} of length $a$ and fundamental domains $\Delta_j$, $1\le j\le a$, we mean
a negative-Schwarzian diffeomorphism 
$$
f:\Delta_1\cup\Delta_2\cup\cdots\cup\Delta_a\,\to\,
\Delta_2\cup\Delta_3\cup\cdots\cup\Delta_{a+1}
$$
such that $f(\Delta_j)=\Delta_{j+1}$. 

The basic geometric estimate on almost parabolic maps is due to J.-C.~Yoccoz.

\proclaim{Yoccoz's Lemma} Let $f:I\to J$ be an almost parabolic map of length
$a$ and fundamental domains $\Delta_j$, $1\le j\le a$. If $|\Delta_1|\ge
\sigma|I|$ and $|\Delta_a|\ge \sigma|I|$, then  
$$
\frac{1}{C_\sigma}\frac{|I|}{\min\{j, a-j\}^2}\;\le\;|\Delta_j|\;\le\;
C_\sigma\frac{|I|}{\min\{j, a-j\}^2}
\ ,
$$
where $C_\sigma>1$ depends only on $\sigma$.
\endproclaim

For a proof, see Appendix B.
We will use Yoccoz's estimates to compare two almost parabolic maps.

\proclaim{Proposition {4.1}} Let $f$ and $g$ be two almost parabolic maps
with the same length $a$ defined on the same interval. Then for all 
$x\in\Delta_1(f)\cap \Delta_1(g)$ and all $0\le k\le a/2$ we have
$$
|f^k(x)-g^k(x)|\;\le\;C\,\|f-g\|_{C^0}k^3
\ .
\slabel{cube}
$$
\endproclaim
\demo{Proof} First note, using the mean-value theorem, that 
$$
\eqalign{
|f^k(x)-g^k(x)|\;&=\;\left|\sum_{j=0}^{k-1}\left(f^{k-j-1}(f(g^j(x))) - 
f^{k-j-1}(g^{j+1}(x))\right)\right|\cr
{}&{}\cr
&\le\;\sum_{j=0}^{k-1}\left|Df^{k-j-1}(\xi_j)\right|\,
\left|f(g^j(x))-g(g^j(x))\right|\ ,\cr}
$$
where $\xi_j$ lies between $f(g^j(x))$ and $g^{j+1}(x)$. Hence we have
$$
|f^k(x)-g^k(x)|\;\le\;\|f-g\|_{C^0}\,
\sum_{j=0}^{k-1}\left|Df^{k-j-1}(\xi_j)\right|
\ .
\slabel{xis}
$$

Let us estimate each summand in the right-hand side of {\eqnoxis}. 
Let $m=m(j)$ be such that $\xi_j\in \Delta_{j+m}$, and assume also that
$j+m\le a/2$. This last condition is always satisfied if the central
fundamental domain of $g$ lies to the left of the central fundamental domain
of $f$ (if this is not the case, then reverse the roles of $f$ and $g$ in
{\eqnoxis} and throughout). Using Yoccoz's Lemma, we see that
$$
|Df^{k-j-1}(\xi_j)|\;\asymp\;\frac{(j+m)^2}{(a-k-m+1)^2}
\;\le\;\left(\frac{j+m}{j+1}\right)^2
\ .
\slabel{fxis}
$$
Hence, it suffices to estimate $m$ as a function of $j$. For this purpose,
let $n=n(j)$ be such that $g^{j+1}(x)\in [f^{j+n-1}(x),f^{j+n}(x)]$. We claim
that $m\le n+1$. There are two possibilities. The first is that $f(g^j(x))\ge
g^{j+1}(x)$: in this case we see easily that 
$$ 
\xi_j\in [g^{j+1}(x),f(g^j(x))]\subseteq [f^{j+n-1}(x),f^{j+n+1}(x)]
$$
and so $m\le n+1$. The second is that $f(g^j(x))<g^{j+1}(x)$. In this case we
have $\xi_j<g^{j+1}(x)<f^{j+n}(x)\in\Delta_{j+n+1}$, so once again $m\le
n+1$. This proves our claim.

$$
\psannotate{
\psboxto(10.4cm;0cm){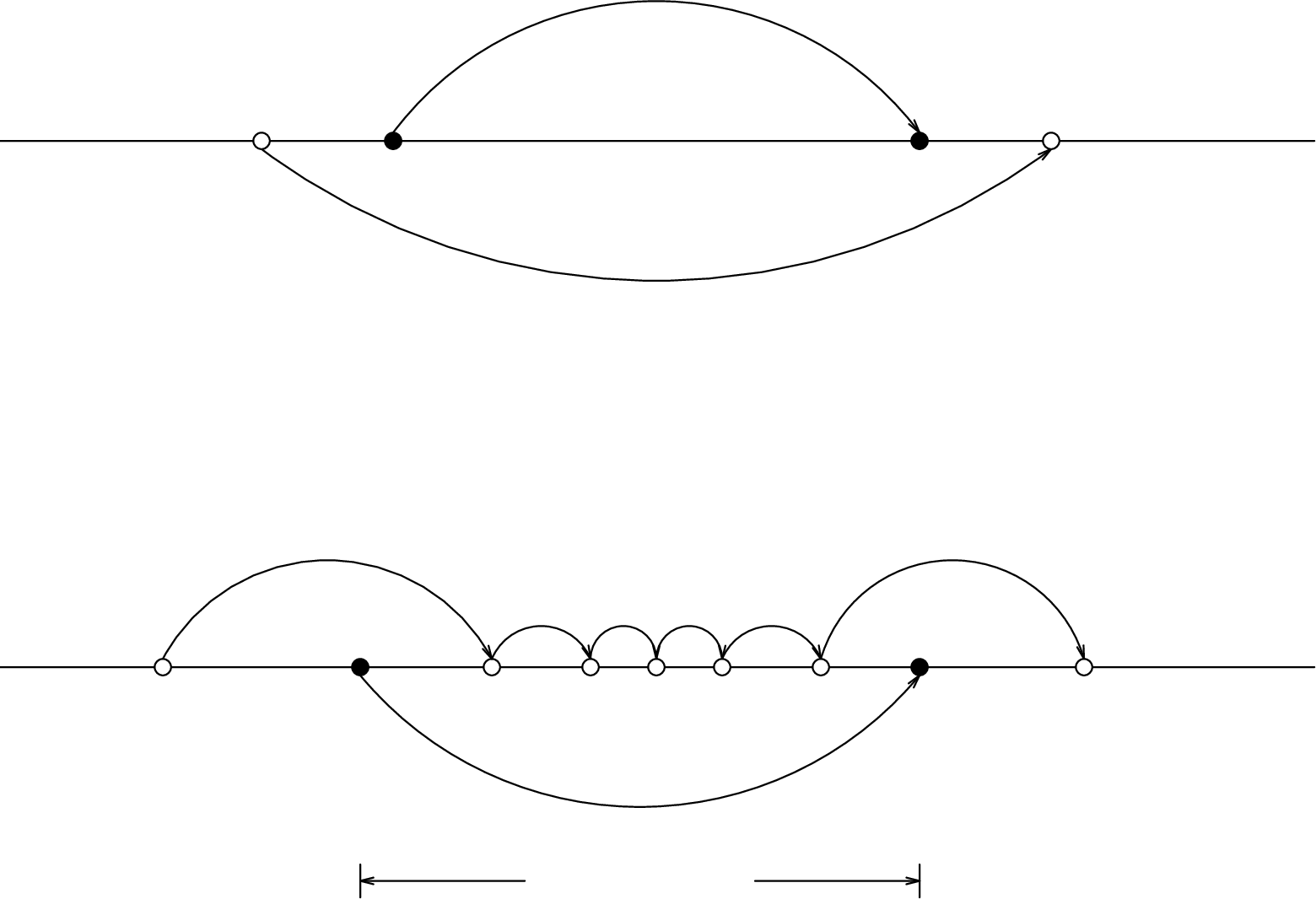}}
{\fillinggrid\at{15.4\pscm}{0.2\pscm}{$\Delta$}
             \at{3.9\pscm}{3.4\pscm}{$f^{j+n-1}(x)$}
             \at{8.35\pscm}{4.85\pscm}{$g^{j+1}(x)$}
             \at{10.7\pscm}{3.5\pscm}{$f^{j+n}(x)$}
             \at{19.2\pscm}{4.85\pscm}{$g^{j+2}(x)$}
             \at{22\pscm}{3.4\pscm}{$f^{j+n+p}(x)$}
             \at{0.5\pscm}{7.5\pscm}{($b$)}
             \at{0.1\pscm}{17.2\pscm}{($a$)}
             \at{3.0\pscm}{13.5\pscm}{$f^{j+n-1}(x)$}
             \at{20.6\pscm}{13.5\pscm}{$f^{j+n}(x)$}
             \at{5.85\pscm}{15.5\pscm}{$g^{j+1}(x)$}
             \at{17.9\pscm}{15.5\pscm}{$g^{j+2}(x)$}
}
$$
\bigskip
\centerline{\sl Figure 3. Bounding $n$ as a function of $j$.}
\bigskip

So now we must bound $n$ as a function of $j$. Again, there are two cases to
consider. 

\itemitem{($a$)} We have $[g^{j+1}(x),g^{j+2}(x)]\subseteq
[f^{j+n-1}(x),f^{j+n}(x)]$. In this case (Figure 3$a$) Yoccoz's Lemma
gives us  
$$
\frac{1}{j^2}\;\le\;\frac{C}{(j+n)^2}
\ ,
$$ 
which implies $n\le Cj$. 

\itemitem{($b$)} We have $g^{j+2}(x)>f^{j+n}(x)$. In this case (Figure
3$b$), $f^{j+n}(x)$ is the first point in the $f$-orbit of $x$ 
that lands inside the interval $\Delta=[g^{j+1}(x),g^{j+2}(x)]$. Let
$p$ be such that $f^{j+n+i}(x)\in \Delta$ for $i=0,1,\ldots,p-1$ but
$f^{j+n+p}(x)\notin \Delta$. Then we have $\Delta\subseteq
[f^{j+n-1}(x), f^{j+n+p}(x)]$, and this time Yoccoz's Lemma gives us 
$$
\frac{1}{j^2}\;\le\;C\,\left(\frac{1}{(j+n)^2}+\frac{1}{(j+n+1)^2}+ \cdots + 
\frac{1}{(j+n+p)^2}\right)\;\le\;\frac{C}{j+n}
\ .
$$
Therefore $n\le Cj^2$ in this case.

In either case we see that $m\le Cj^2$. Carrying this information back to
{\eqnofxis}, we deduce that 
$$
|Df^{k-j-1}(\xi_j)|\;\le\;C\,j^2
\ .
\slabel{fsquare}
$$
Substituting {\eqnofsquare} into {\eqnoxis}, we arrive at {\eqnocube}, and the
proof is complete.\qed
\enddemo

\subhead
4.2 A criterion for smoothness
\endsubhead
One key ingredient in the proof of our First Main Theorem is a slight
extension of a result originally due to Carleson {\refCa}, namely Proposition
{4.3} below.
To formulate it, we need an auxiliary definition.

\demo{Definition} A {\it fine grid} is a sequence $\{{\Cal Q}_n\}_{n\ge 0}$ of
finite partitions of $S^1$ which satisfies
\itemitem{($a$)} Each ${\Cal Q}_{n+1}$ is a
strict refinement of  ${\Cal Q}_n$; 
\itemitem{($b$)} There exists $a>0$ such that each $I\in {\Cal Q}_n$ is the
disjoint union of at most $a$ atoms of ${\Cal Q}_{n+1}$; 
\itemitem{($c$)} There exists $c>0$ such that $c^{-1}|I|\le |J|\le c|I|$ for
each pair of adjacent atoms $I,J\in {\Cal Q}_n$. 
\enddemo

For example, the dynamical partitions $\{{\Cal P}_n\}$ of a critical circle
map with rotation number of bounded type always form a fine grid, by Theorem
{3.1}. We note the following easy lemma concerning a fine grid $\{{\Cal
Q}_n\}$.  

\proclaim{Lemma {4.2}} If $I\in {\Cal Q}_n$, $J\in {\Cal Q}_{n+1}$
and $J\subseteq I$, then $(1+c^{-1})|J|\le |I|\le ac^a|J|$. In
particular, there exist $C_0>1$ and $0<\lambda_0<\lambda_1<1$ such that
$C_0^{-1}\lambda_0^n\le |I|\le C_0\lambda_1^n$, for all $I\in {\Cal
Q}_n$. \qed 
\endproclaim

The constants $a$, $c, C_0, \lambda_0, \lambda_1$, are the {\it fine
constants} of $\{{\Cal Q}_n\}$.

\proclaim{Proposition {4.3}} Let $h:S^1\to S^1$ be a homeomorphism and let
$\{{\Cal Q}_n\}_{n\ge 0}$ be a fine grid. 
\itemitem{($a$)} If there exists $C>0$ such that
$$
\left|\frac{|I|}{|J|}-\frac{|h(I)|}{|h(J)|}\right|\;\le\; C
\ ,
$$
for each pair of adjacent atoms $I,J\in {\Cal Q}_n$, for all $n\ge 0$,
then $h$ is quasisymmetric. 
\itemitem{($b$)} If there exist constants $C>0$ and $0<\lambda <1$ such that  
$$
\left|\frac{|I|}{|J|}-\frac{|h(I)|}{|h(J)|}\right|\;\le\; C\lambda^n
\ ,
\slabel{mesh}
$$
for each pair of adjacent atoms $I,J\in {\Cal Q}_n$, for all
$n\ge 0$, then $h$ is a $C^{1+\alpha}$-diffeomor\-phism for some
$\alpha>0$. 
\endproclaim

The proof of Proposition {4.3} will depend on the following fact from
elementary real analysis. If $\phi$ is a real-valued function
in an interval or oriented arc on the circle, let
$D^+\phi(x)=\lim_{t\downarrow 0}{(\phi(x+t)-\phi(t))/t}$ be the {\it
right derivative} of $\phi$ at $x$, if the limit exists.

\proclaim{Lemma {4.4}} Let $\phi_n:[0,1]\to {\Bbb R}$ be a sequence of
continuous, right differentiable mappings such that $D^+\phi_n$
converges uniformly to an $\alpha$-H\"older continuous function
$\varphi:[0,1]\to {\Bbb R}$, and such that each $D^+\phi_n$ is
Riemann-integrable. If $\phi_n$ converges uniformly to $\phi$, then
$\phi$ is $C^{1+\alpha}$ and $D\phi=\varphi$. \qed
\endproclaim

\demo{Proof of Proposition {4.3}} 
We will prove ($b$) only, the proof of ($a$) being somewhat easier.
Let $\phi_n$ be the piecewise affine
$C^0$-approximations to $h$ determined by the vertices of ${\Cal Q}_n$.
Then $\phi_n$ is differentiable on the right, and $D^+\phi_n$ is a
step function. First we show that $\left\{D^+\phi_n\right\}_{n\ge 0}$
is a uniform Cauchy sequence, and then that the limit is H\"older 
continuous. 
Take an atom $I$ of ${\Cal Q}_n$, and consider the decomposition
$$
I\;=\;J_1\cup J_2\cup\cdots \cup J_p
\ ,
$$
with $J_k\in {\Cal Q}_{n+1}$ consecutive and pairwise disjoint  and
$p\le a$. Then $D^+\phi_n$ is constant on $I$ and $D^+\phi_{n+1}$ is
constant on each $J_k$, say
$$
\cases
D^+\phi_n(t)\;=\;\sigma\;=\;\frac{\d |\phi_n(I)|}{\d |I|} &{\ (t\in I)}\cr
{}&{}\cr
D^+\phi_{n+1}(t)\;=\;\sigma_k\;=\;\frac{\d |\phi_{n+1}(J_k)|}{\d |J_k|} &{\
(t\in  J_k)}\cr  
\endcases
\ .
$$
Thus, we have
$$
\sigma|I|\;=\;\sum_{k=1}^p \sigma_k|J_k|
\ ,
$$
and in particular $\sigma'=\min \sigma_k\le \sigma\le \max
\sigma_k=\sigma''$. Also, $\sigma'/\sigma''\le
\sigma/\sigma_k\le \sigma''/\sigma'$ for all $k$. Since by assumption
$|1-(\sigma_{k+1}/\sigma_k)| \le C\lambda^{n+1}$, an easy
telescoping trick gives us 
$$
\frac{\sigma''}{\sigma'}\;\le\;(1+C\lambda^{n+1})^{a}\;\le\;1+C\lambda^{n+1}
\ .
$$
A similar lower bound holds true for $\sigma'/\sigma''$. Therefore we have 
$$
1-C\lambda^n\;\le\;\frac{\sigma}{\sigma_k}\;\le\;1+C\lambda^n
\ ,
\slabel{slope}
$$
for all $k=1,2,\ldots p$. This shows that the sequence
$\left\{D^+\phi_n\right\}_{n\ge 0}$ is uniformly bounded, and moreover
that for all $m\ge n\ge 0$ and all $t\in S^1$, we have
$$
\left|D^+\phi_m(t)-D^+\phi_n(t)\right|\;\le\;
C\sum_{j=n}^{m-1}\lambda^j\;<\;\frac{C}{1-\lambda}\,\lambda^n
\ .
\slabel{cauchy}
$$
Hence $\left\{D^+\phi_n\right\}_{n\ge 0}$ is a uniform Cauchy
sequence as claimed. Let $\varphi=\lim{D^+\phi_n}$, and let $\alpha>0$
be such that $\lambda_0^{\alpha}=\lambda$. We prove $\varphi$ is
$\alpha$-H\"older as follows.
It suffices to consider points $x,y\in S^1$ whose distance is smaller
than $\inf_{I\in {\Cal Q}_0}|I|$. Take the smallest $n$ such that $x$
and $y$ belong to distinct elements of ${\Cal Q}_n$. Then either $n=0$ or $x$
and $y$ lie in a common element of ${\Cal Q}_{n-1}$. Either way we have by
{\eqnoslope} 
$$
\left|D^+\phi_n(x)-D^+\phi_n(y)\right|\;\le\;C\lambda^n
\ .
\slabel{holder}
$$
Combining {\eqnocauchy} and {\eqnoholder}, we deduce that 
$$
\eqalign{
\left|\varphi(x)-\varphi(y)\right|&\le\left|\varphi(x)-D^+\phi_n(x)\right| 
+\left|D^+\phi_n(x)-D^+\phi_n(y)\right|+\left|D^+\phi_n(y)-\varphi(y)\right|\cr
&{}\cr
&\le\;\frac{C}{1-\lambda}\lambda^n+C\lambda^n+\frac{C}{1-\lambda}\lambda^n
\;\le\;C\lambda_0^{n\alpha}\cr
&{}\cr
&\le\;C|x-y|^{\alpha}\ ,\cr}
$$
and so $\varphi$ is $\alpha$-H\"older as claimed.\qed
\enddemo

\demo{Remark} In the language of conditional expectations, the
sequence $\{D^+\phi_n\}_{n\ge 0}$ satisfies 
${\Bbb E}\left(D^+\phi_n\,|\,{\Cal B}_n\right)=D^+\phi_{n+1}$, where
${\Cal B}_n$ is the $\sigma$-algebra generated by ${\Cal Q}_n$, and
therefore constitutes a {\it martingale}. Thus, the existence of a
pointwise $a.e.$ limit $\varphi$, merely as an integrable function, is
a special case of J.~Doob's {\it martingale convergence theorem}, see
{\refBi}, p.~490.  
\enddemo

\subhead
4.3 A suitable fine grid
\endsubhead
The dynamical partitions ${\Cal P}_n$ of a critical circle map $f$ do
{\it not} determine a fine grid, unless the rotation number of $f$ is
of bounded type. We will however use these dynamical partitions
to build a fine grid $\{{\Cal Q}_n\}$ for our map $f$.
The construction requires some preliminary definitions.

An element $I\in {\Cal P}_n$ is a {\it saddle-node}
atom if it is the disjoint union of some number $a\ge 1000$ of atoms
of ${\Cal P}_{n+1}$. 

Given two atoms ${\Cal P}_{n+1}\ni J\subseteq I\in
{\Cal P}_n$, the {\it order} of $J$ inside $I$ is one plus the
smallest number of atoms of ${\Cal P}_{n+1}$ on the right and left
components of $I\setminus J$. 

Note that inside a saddle-node atom $I\in {\Cal P}_n$ there are exactly
two atoms of ${\Cal P}_{n+1}$ of order $k$ for each $k\le a/2$. Let $N\ge 0$
be such that $2^{N+1}<a/2$. For each $0\le i\le N$, we define $M_i$, the {\it
$i$-th central interval} of $I$, to be the convex-hull
$[J,J^*]\subseteq I$ of the union of both atoms $J,J^*$ of order
$2^i$. Note that these central intervals are nested (see Figure 4).
The left and right components of $M_i\setminus M_{i+1}$, respectively
$L_i$ and $R_i$, are the {\it lateral intervals} of $I$. The central
interval $M_N$ is also called the {\it final interval} of
$I$. The lateral intervals together with the final interval form
a special partition of $I$, the {\it balanced} partition of $I$. 

\demo{Remark} It follows from Yoccoz's lemma that $|L_i|\asymp
|M_{i+1}|\asymp |R_i|$ for all $i$.
\enddemo

$$
\psannotate{
\psboxto(10cm;0cm){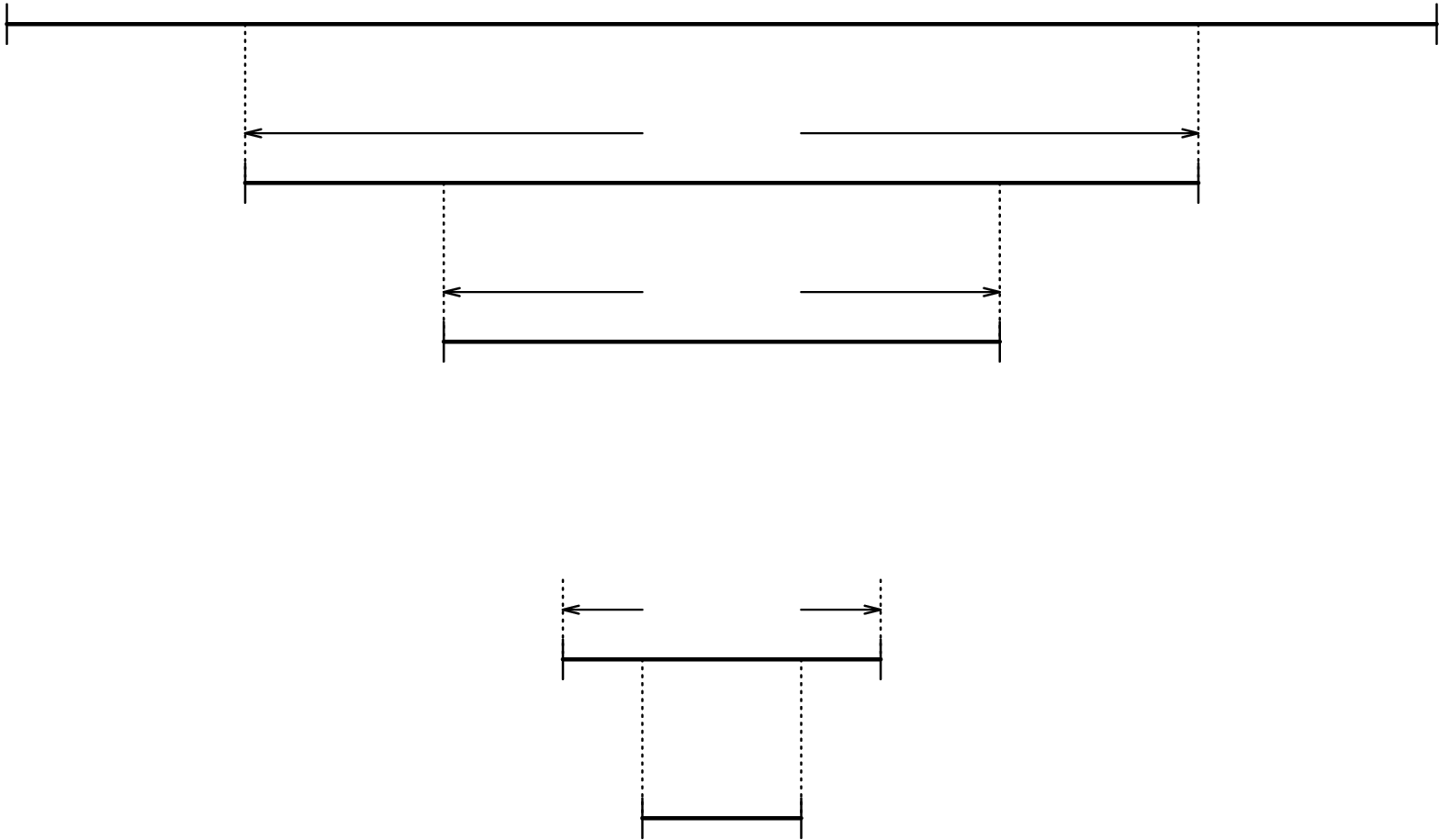}}
{\fillinggrid\at{13.8\pscm}{-0.8\pscm}{$M_{N}$}
             \at{13\pscm}{3.5\pscm}{$M_{N-1}$}
             \at{13.3\pscm}{8.5\pscm}{$M_{2}$}
             \at{13.0\pscm}{11\pscm}{$M_{1}$}
             \at{13.22\pscm}{5.2\pscm}{$\scriptstyle\bullet$}
             \at{12.95\pscm}{6\pscm}{$\scriptstyle\bullet$}
             \at{12.7\pscm}{6.8\pscm}{$\scriptstyle\bullet$}
             \at{2.4\pscm}{12\pscm}{$L_0$}
             \at{21.3\pscm}{12\pscm}{$R_0$}
             \at{5.4\pscm}{9.2\pscm}{$L_1$}
             \at{17.4\pscm}{9.2\pscm}{$R_1$}
}
$$
\bigskip
\centerline{\sl Figure 4. Central and lateral intervals of a saddle-node
             atom.} 
\bigskip
Now we define an auxiliary partition $\widetilde{\Cal P}_n$, for each
$n\ge 1$. The atoms of $\widetilde{{\Cal P}_n}$ are all atoms of ${\Cal
P}_n$ which are not saddle-node, together with the atoms of the
balanced partitions of all saddle-node atoms of ${\Cal P}_n$. The
partition ${\Cal Q}_n$ that we want is constructed from $\widetilde{\Cal
P}_n$ and ${\Cal P}_n$ as follows.

\proclaim{Proposition {4.5}} There exists a fine grid $\{{\Cal Q}_n\}$ in
$S^1$ with the following properties.
\item{($a$)} Every atom of ${\Cal Q}_n$ is the union of at most 3
atoms of ${\Cal Q}_{n+1}$.
\item{($b$)} Every atom $\Delta$ of ${\Cal Q}_n$ is a union of atoms of
${\Cal P}_m$ for some $m\le n$, and there are four possibilities:
\itemitem{($b_1$)} $\Delta$ is a single atom of ${\Cal P}_m$;
\itemitem{($b_2$)} $\Delta$ is a central interval of $\widetilde{\Cal
P}_m$; 
\itemitem{($b_3$)} $\Delta$ is the union of at least two atoms of
${\Cal P}_{m+1}$ contained in a single atom of $\widetilde{\Cal P}_m$.
\itemitem{($b_4$)} $\Delta$ is a union of intervals which are
simultaneously atoms of ${\Cal P}_m$ and $\widetilde{\Cal P}_m$.
\endproclaim

\demo{Proof} The proof is by induction on $n$. The first partition ${\Cal
Q}_1$ consists of all atoms of ${\Cal P}_1$ which are not saddle-node atoms
together with the intervals $L_0$, $M_1$ and $R_0$ of each saddle-node
interval $I\in {\Cal P}_1$ ($I=L_0\cup M_1\cup R_0$). It is clear that each
atom of ${\Cal Q}_1$ falls within one of the categories ($b_1$)-($b_4$)
above. 

Assuming ${\Cal Q}_n$ defined, define ${\Cal Q}_{n+1}$ as follows. Take an
atom $I\in {\Cal Q}_n$ and consider the four cases below.

\item{(1)} If $I$ is a single atom of ${\Cal P}_m$ then one of two things can
happen:
\itemitem{(i)} $I$ is a saddle-node atom: In this case write $I=L_0\cup
M_1\cup R_0$ as above and take $L_0$, $R_0$ and $M_1$ as atoms of ${\Cal
Q}_{n+1}$. Note that the lateral intervals $L_0$ and $R_0$ are atoms of type
($b_1$), while the central interval $M_1$ is of type ($b_2$).
\itemitem{(ii)} $I$ is not a saddle-node atom: In this case write $I=L\cup
M\cup R$ where $L$ and $R$ are the atoms of ${\Cal P}_{m+1}$ adjacent to the
endpoints of $I$ and $M$ is the union of the other atoms of ${\Cal P}_{m+1}$
inside $I$. Add these three intervals to ${\Cal Q}_{n+1}$, noting that $L$ and
$R$ are of type ($b_1$), while $M$ is of type ($b_4$).

\item{(2)} If $I$ is a central interval of $\widetilde{\Cal P}_m$ which is not the
final interval, consider the next central interval of $\widetilde{\Cal P}_m$
inside $I$, say $M$, and the two corresponding lateral intervals $L$ and $R$
such that $I=L\cup M\cup R$, and declare $L$, $R$ and $M$ members of ${\Cal
Q}_{n+1}$. Note that $L$ and $R$ are of type ($b_3$), while $M$ is of type
($b_2$).

\item{(3)} If $I$ is a union of $p\ge 2$ consecutive atoms
$\Delta_1,\ldots,\Delta_p$ of $\widetilde{\Cal P}_{m+1}$ inside a single atom of
${\Cal P}_m$, divide it up into three approximately equal parts. More
precisely, write $p=3q+r$ and, when $r=0$ or $1$, consider $I=L\cup M\cup R$
where 
$$
L\;=\;\bigcup_{j=1}^q\Delta_j \ ,\ M\;=\;\bigcup_{j=q+1}^{p-q}\Delta_j \ ,
\ R\;=\;\bigcup_{j=p-q+1}^p\Delta_j
\ .
$$
When $r=2$, consider $I=L\cup M\cup R$ where
$$
L\;=\;\bigcup_{j=1}^{q+1}\Delta_j \ ,\ M\;=\;\bigcup_{j=q+2}^{p-q-1}\Delta_j
\ ,\ R\;=\;\bigcup_{j=p-q}^p\Delta_j
\ .
$$
Note that $M$ is empty when $p=2$. In any case, we obtain two or three new
atoms of ${\Cal Q}_{n+1}$ which are either single atoms of ${\Cal P}_{m+1}$,
and therefore of type ($b_1$), or once again intervals of type ($b_3$).

\item{(4)} If $I$ is a union of intervals which are simultaneously atoms of
${\Cal P}_m$ and $\widetilde{\Cal P}_m$, divide it up exactly as in (3),
obtaining either two or three new atoms of ${\Cal Q}_{n+1}$ which are either
single atoms of ${\Cal P}_m$, and therefore of type ($b_1$), or once again
intervals of type ($b_4$).

This completes the induction. That $\{{\Cal Q}_n\}_{n\ge 0}$ constitutes a
fine grid follows easily from the real bounds and the remark preceding this
proposition. \qed 

\enddemo

An immediate consequence of the mere existence of such a fine grid is the fact
that any two critical circle maps with the same rotation number are
quasisymme\-tri\-cally conjugate.

\proclaim{Corollary {4.6}} Let $f$ and $g$ be critical circle maps with
the same irrational rotation number, and let $h$ be the conjugacy between $f$
and $g$ that maps the critical point of $f$ to the critical point of $g$. Then
$h$ is quasisymmetric.
\endproclaim
\demo{Proof} Apply Proposition {4.3} ($a$) to the fine grid
constructed above. \qed 
\enddemo

\subhead
4.4 Proof of the First Main Theorem
\endsubhead
The strategy for the proof of our First Main Theorem is as follows. Given
two critical circle maps $f$ and $g$ with the same rotation number, consider
the special partitions ${\Cal Q}_n={\Cal Q}_n(f)$ and $\widetilde{\Cal Q}_n={\Cal
Q}_n(g)$ given by Proposition {4.5}. The conjugacy $h$ is an isomorphism
between the corresponding fine grids, and we want to show that the coherence
property {\eqnomesh} holds for $h$ and $\{{\Cal Q}_n\}$. To do this, we use
the fact that the renormalizations $f_n$ and $g_n$ are exponentially close to
prove that the vertices of ${\Cal Q}_{n+p}$ are exponentially close to the
corresponding vertices of $\widetilde{\Cal Q}_{n+p}$, provided $p$ is a small
fraction of $n$. This property is proved in two steps: first for those
vertices of ${\Cal Q}_n$ that lie in the domain of $f_n$, a step that works
for arbitrary rotation numbers, and then propagated via Koebe's distortion
principle to the other vertices. The first step holds without restriction on
the rotation number, but the second does not. And it could not, indeed, as the
counterexamples of section 5 will show.

In what follows we use the notation $x_n=x_n(f)=f^{q_n}(x)$. If the
renormalizations $f_n$ and $g_n$ converge together exponentially fast then
$x_n(f)/x_n(g)\to 1$ and $|x_n(f)-x_n(g)|/|x_n(f)|\to 0$ exponentially fast
as well.

\demo{Definition} Let $f_m:J_m\to J_m$ be the $m$-th renormalization of $f$
and let $k$ be an integer such that $|k|\le a_m/2$. The {\it restricted
domain} of $f_m^k$, denoted $D_{m,k}$, is defined as follows.
$$
D_{m,k}\;=\;
\cases
I_m\cup [f^{\frac{a_m}{2}-k}(x_{m-1}), x_{m-1}],&\text{when $k>0$}\cr
{}&{}\cr
\left(I_m\setminus [x_m,f_m(x_{m+1})]\right)\cup
[c,f_m^{\frac{a_m}{2}}(x_{m-1}) ],&\text{when $k=-1$}\cr
{}&{}\cr
\left(I_m\setminus [x_m,f_m(x_{m+1})]\right)\cup
[c,f_m^{-\frac{a_m}{2}-k}(x_{m+1}) ],&\text{when $k<-1$}.\cr
\endcases
$$
\enddemo

\proclaim{Lemma {4.7}} For all $x\in D_{m,k}$ we have $|Df_m^k(x)|\le C$,
where $C>0$ depends only on the real bounds.
\endproclaim
\demo{Proof} Use Theorem {3.1} and Yoccoz's Lemma.\qed
\enddemo

\proclaim{Lemma {4.8}} Let $v_0$ be a vertex of ${\Cal P}_{n+p}$ such
that $v_0\in J_n$. Then there exist $n\le m\le n+p$ and $1\le N\le p$ such
that $v_0$ can be represented in the form
$$
v_0\;=\;\phi_1\circ\phi_2\circ\cdots\circ\phi_N(x_m)
\ ,
$$
where $\phi_j=f_{m_j}^{k_j}$ for some $n\le m_j\le n+p$ and $|k_j|\le
a_{m_j}/2$, and where the point $\phi_{j+1}\circ\cdots\circ\phi_N(x_m)$
belongs to the restricted domain of $\phi_j$ for each $j$. Moreover, if $v_0$
is also a vertex of ${\Cal Q}_{n+p}$, then $|k_j|\le 2^p$ for all $j$. 
\endproclaim

\demo{Proof} Let $n\le m_1\le n+p$ be largest with the property that $v_0\in
J_{m_1}\setminus J_{m_1+1}$, and let $0<i\le a_{m_1}$ be such that
$f_{m_1}^i(v_0)\in J_{m_1+1}$. If $i\le a_{m_1}/2$ then let $k_1=-i$;
otherwise let $k_1=i-a_{m_1}$. We get $\phi_1=f_{m_1}^{k_1}$ and a new vertex
$v_1=f_{m_1}^{-k_1}(v_0)\in J_{m_1+1}$. If $v_1\in J_{n+p}$ then
$v_1=f_{n+p}(x_{n+p-1})$ necessarily, and we can stop. On the other hand, if
$v_1\notin J_{n+p}$, then once again there exists $m_2$ in the range $m_1<
m_2< n+p$ such that $v_1\in J_{m_2}\setminus J_{m_2+1}$, and we can proceed
inductively. At the end of this process we get sequences
$m_1<m_2<\cdots<m_N\le n+p$ (so $N\le p$) and $v_1, v_2,\ldots,v_N$ with
$v_j\in J_{m_j}\setminus J_{m_j+1}$, and for each $j$ an integer $k_j$ with
$|k_j|\le a_{m_j}/2$ such that $v_{j+1}=f_{m_j}^{-k_j}(v_j)$. The last vertex
$v_N$ is necessarily $x_m$ for some $m\le n+p$. Hence it suffices to take
$\phi_j=f_{m_j}^{k_j}$ to get the desired representation. The last statement
of the Lemma is immediate from the way the partitions ${\Cal Q}_n$ were
constructed. \qed
\enddemo

Now we want to estimate $|v_0-w_0|$, where $v_0\in J_n(f)$ is a vertex of
${\Cal Q}_{n+p}(f)$ and $w_0$ is the corresponding vertex of ${\Cal
Q}_{n+p}(g)$. Here we assume $p\le \sigma n$ for some small $\sigma>0$ to be 
chosen later. By Lemma {4.8} above, there exist points $x_m=x_m(f)$,
$y_m=x_m(g)$ and a number $N\le p$ such that 
$$
|v_0-w_0|\;=\;\left|\phi_1\circ\phi_2\circ\cdots\circ\phi_N(x_m) -
\psi_1\circ\psi_2\circ\cdots\circ\psi_N(y_m)\right|
\ ,
$$ 
where $\phi_j=f_{m_j}^{k_j}$ and $\psi_j=g_{m_j}^{k_j}$, with $n\le m_j\le
n+p$ and $|k_j|\le 2^p$. 

For all $n$, the return map $f_n$ is an almost parabolic map on
$J_{n-1}\setminus J_n$, and similarly for $g_n$. Our hypothesis is that
$\|f_n-g_n\|_0 \le C\lambda^n$ for all $n$, for some $0<\lambda<1$. If $x$ is a
point in the first fundamental domain of both $f_{m_j}$ and $g_{m_j}$, then
by Proposition {4.1} we have
$$
\eqalign{
\left|\phi_j(x)-\psi_j(x)\right|\;&\le\;C\lambda^{m_j}k_j^3
\;\le\;C\lambda^n2^{3p}\cr 
{}&{}\cr
&\le\;C\left(\lambda 2^{3\sigma}\right)^n\;=\;C\theta^n\ ,}
$$
where $\theta=\lambda 2^{3\sigma}<1$ if $\sigma$ is small enough.

Using this, and since $|x_m-y_m|\le C\lambda^n$, we see that
$$
\eqalign{
|\phi_N(x_m)-\psi_N(y_m)|\;&\le\;
|\phi_N(x_m)-\psi_N(x_m)|+|\psi_N(x_m)-\psi_N(y_m)|\cr
&{}\cr
&\le\;  
C\theta^n + \|D\psi_N\|_0\,|x_m-y_m|
\;\le\; C\left(\theta^n+\lambda^n\right)\;\le\; C\theta^n\ ,}
$$
because $\|D\psi_N\|_0$ is bounded, by Lemma {4.7}. Proceeding
inductively, we get 
$$
\eqalign{
\left|\phi_1\circ\cdots\circ\phi_N(x_m)-
\psi_1\circ\cdots\circ\psi_N(y_m)\right|\;&\le\;
\left(C+C^2+\cdots+C^N\right)\theta^n\cr 
{}&{}\cr
&\le\;C^p\theta^n\le (C^\sigma\theta)^n\ ,}
$$
so making $\sigma$ still smaller, $|v_0-w_0|$ is exponentially small as
desired.

\heading 5. 
Counterexamples to $C^{1+\alpha}$ rigidity
\endheading
Our purpose now is to construct $C^\infty$ counterexamples to the conjectured
$C^{1+\alpha}$ rigidity of critical circle maps. We will consider critical
circle maps whose rotation number $\rho(f)=[a_0,a_1,\ldots,a_n,\ldots]$
satisfies
$$
\cases
\lim\sup{\displaystyle\frac{1}{n}}\log{a_n}\;=\;\infty &{}\cr
{}&{}\cr
a_n\;\ge\;2 & \text{for all $n$},\cr
\endcases
\slabel{rot}
$$
The class of all rotation numbers satisfying {\eqnorot} will be denoted by
$S$. It can be shown that the Hausdorff dimension of $S$ is less than or equal
to $1/2$, see {\refGd}. On the other hand, $S$ contains Diophantine numbers:
for example, the number $\rho$ whose partial quotients are $a_n=2^{2^n}$ is
Diophantine and satisfies {\eqnorot}.

\proclaim{Theorem {5.1}} For every $\rho\in S$ there exist $C^\infty$
critical circle maps $f$, $g$ with $\rho(f)=\rho(g)=\rho$ such that $f$ and
$g$ are not $C^{1+\beta}$ conjugate for any $\beta>0$. 
\endproclaim

The proof will make use of a $C^\infty$ surgery procedure that we explain
next. These counterexamples have one additional feature: their
successive renormalizations do converge together at an exponential
rate. This will be clear from the construction.

\subhead
5.1 Saddle-node surgery
\endsubhead
Given $f$ as above and a fixed $n\ge 1$, let $J_n=J_n(f)=[f^{q_n}(c),
f^{q_{n-1}}(c)]\subseteq S^1$ be the $n$-th renormalization interval of
$f$. When $a_n$ is very large, the first return map $f_n:J_n\to J_n$ is an
almost parabolic map of length $a_n$. 

Let $\Delta_1$ be the fundamental domain of this almost parabolic map which
is adjacent to $x_{n-1}=f^{q_{n-1}}(c)$, and let
$\Delta_j=f_n^{j-1}(\Delta_1)$, $j\le a_n$. Let $z\in \Delta_1$ be the point
such that $f_n^{a_n}(z)=x_{n+2}=f^{q_{n+2}}(c)$, that is,
$z=f^{q_{n+2}-a_nq_n}(c)$. Note that since $a_n\ge 2$, $x_{n+2}$ is not an
endpoint of $f_n^{a_n}(\Delta_1)$, and so by the real bounds it splits
$f_n^{a_n}(\Delta_1)$ into two intervals of comparable lengths. Hence the same
holds for $z$. Namely, $z$ splits $\Delta_1$ into two intervals $L$, $R$ with
$|L|\asymp |R|$. In particular we have $\tau |\Delta_1|\le |L|\le
(1-\tau)|\Delta_1|$ (and similarly for $R$) for some constant $\tau$ depending
on the real bounds; we use this fact in the proof of Proposition {5.2} below.

Consider now another critical circle map $\tilde{f}$ with the same rotation
number as $f$, the interval $\tilde{J_n}=J_n(\tilde{f})$, the first return map
$\tilde{f}_n:\tilde{J_n}\to \tilde{J_n}$, the point
$\tilde{z}=\tilde{f}^{q_{n+2}-a_nq_n}(\tilde{c})$ and the corresponding
intervals $\tilde{L}$, $\tilde{R}$. Also, let $N=\lceil a_n/2\rceil$.

\demo{Definition} The number 
$$
\left|\frac{|f_n^{N-1}(L)|}{|f_n^{N-1}(R)|} -
\frac{|\tilde{f}_n^{N-1}(\tilde{L})|}{|\tilde{f}_n^{N-1}(\tilde{R})|}\right|
$$
is called the $n$-th order discrepancy between $f$ and $\tilde{f}$.
\enddemo

\proclaim{Proposition {5.2}} Given a $C^\infty$ critical circle map $f$
with $\rho(f)\in S$, consider a function $\sigma(n)\to\infty$ such that
$$
\lim\sup{\displaystyle\frac{1}{n\sigma(n)}}\log{a_n}\;=\;\infty 
\ .
$$
Then for all $n\ge 1$, there exists a critical circle map
$\tilde{f}=F(n;f)$  with the same rotation number and critical point as $f$
and having the following properties.
\itemitem{(a)} We have $\tilde{f}^j(c)=f^j(c)$ for $0\le j\le q_{n+1}$; in
particular, $J_n(\tilde{f})=J_n=J_n(f)$.
\itemitem{(b)} We have $\tilde{f}=\Phi\circ f$, where $\Phi$ is a $C^\infty$
diffeo such that
$$
\|\Phi^{\pm 1}-\roman{id}_{S^1}\|_{C^k}\le B_k|J_n|^{\sigma(n)-k+1}
$$
for all $k$, where $B_k>0$ is constant depending only on $k$.
\itemitem{(c)} The $n$-th order discrepancy between $f$ and $\tilde{f}$ is
$\ge C|J_n|^{2\sigma(n)}$.
\itemitem{(d)} We have $J_{n+1}(\tilde{f})=J_{n+1}(f)$ and
$\tilde{f}_{n+1}=f_{n+1}$; in particular, $m$-th order discrepancy between $f$
and $\tilde{f}$ is equal to zero for all $m>n$.
\endproclaim

\demo{Proof} We modify $f$ inside $f^{-1}(\Delta_1)$ using a $C^\infty$ bump
function so as to move $z$ by a distance $\ge C|\Delta_1|^{1+\sigma(n)}$ inside
$\Delta_1$. This we do as follows.

Let $\varphi: [0,1]\to [0,1]$ be a $C^\infty$ perturbation of the identity
such that $|\varphi(x)-x|\ge |\Delta_1|^{\sigma(n)}$ for all $\tau \le x\le
1-\tau$ (and $\tau$ as above), and such that $|D^k\varphi(x)|\le
B_k|\Delta_1|^{\sigma(n)}$ for all $0\le x\le 1$ and all $k\ge 2$. Define
$\phi:\Delta_1\to \Delta_1$ by $\phi=A\circ\varphi\circ A^{-1}$ where $A$ is
the affine orientation-preserving map that carries $[0,1]$ onto
$\Delta_1$. Note that $|\phi(z)-z|\ge |\Delta_1|^{1+\sigma(n)}$. Moreover,
since $D^k\phi=|\Delta_1|^{1-k}D^k\varphi$, we have 
$$
\|\phi^{\pm
1}-\roman{id}_{\Delta_1}\|_{C^k}\;\le\;B_k|\Delta_1|^{\sigma(n)-k+1} 
$$
for all $k$. Define $\psi:\Delta_{a_n}\to \Delta_{a_n}$ as the conjugate of
$\phi^{-1}$ by the diffeo $f_n^{a_n-1}: \Delta_1\to \Delta_{a_n}$, namely
$$
\psi\;=\; f_n^{a_n-1}\circ \phi^{-1}\circ (f_n^{a_n-1})^{-1}
\ .
\slabel{psiphi}
$$
Using the $C^m$ Approximation Lemma (see Appendix A), we see from {\eqnopsiphi}
that  
$$
\|\psi^{\pm 1}-\roman{id}_{\Delta_{a_n}}\|_{C^{k-1}}
\;\le\; C\|\phi^{\pm 1}-\roman{id}_{\Delta_1}\|_{C^k}
\;\le\; B_k|\Delta_1|^{\sigma(n)-k+1}
\ .
$$
Define $\Phi:S^1\to S^1$ to be equal to $\phi$ on $\Delta_1$, to $\psi$ on
$\Delta_{a_n}$ and to the identity everywhere else. The critical circle map we
look for is $\tilde{f}=\Phi\circ f$. 
Note that $\|\Phi^{\pm 1}-\roman{id}_{S^1}\|_{C^k}\le
B_k|\Delta_1|^{\sigma(n)-k+1}$ for all
$k$; since $|\Delta_1|\asymp |J_n|$ by the real bounds, this proves ($b$).
It is also clear from the construction that property ($a$) holds too. It
follows in particular that the first $n+1$ partial quotients of the rotation
number of $\tilde{f}$ agree with those of $f$. More remarkable is that, because
what $\phi$ does is undone by $\psi$, we have
$$
\cases
\tilde{f}^{q_n}|I_{n+1}=f^{q_n}|I_{n+1}&{}\cr
{}&{}\cr
\tilde{f}^{q_{n+1}}|I_n=f^{q_{n+1}}|I_n&{}\cr
\endcases
\ .
$$
In other terms, $\tilde{f}_n=f_n$, the $n$-th renormalizations
agree. Therefore all subsequent renormalizations agree as well. This shows
that $\rho(\tilde{f})=\rho(f)$ and also proves ($d$).

It remains to prove ($c$), so we estimate the $n$-th order discrepancy between
$f$ and $\tilde{f}$ from below. Since $|z-\tilde{z}|\ge
|\Delta_1|^{1+\sigma(n)}$, a 
simple calculation yields 
$$
\left|\frac{|L|}{|R|}-\frac{|\tilde{L}|}{|\tilde{R}|}\right|\;\ge\;
C|\Delta_1|^{\sigma(n)}\;\ge\; C|J_n|^{2\sigma(n)}
\ ,
\slabel{discrep}
$$
provided $n$ is sufficiently large. Since, by the real bounds, the map
$f_n^{N-1}:\Delta_1\to \Delta_N$ has bounded distortion, and since
$\tilde{f}_n=f_n$, {\eqnodiscrep} gives us 
$$
\left|\frac{|f_n^{N-1}(L)|}{|f_n^{N-1}(R)|} -
\frac{|\tilde{f}_n^{N-1}(\tilde{L})|}{|\tilde{f}_n^{N-1}(\tilde{R})|}\right|
\;\ge\; C|J_n|^{2\sigma(n)}
\ ,
$$
and this proves ($c$).\qed
\enddemo

\subhead 
5.2 The counterexamples
\endsubhead
We now iterate the procedure given by Proposition {5.2} to prove our Second
Main Theorem (that is, Theorem {5.1}). We start with a $C^\infty$
map $f$ with $\rho(f)\in S$ as before and select $n_1<n_2 <\cdots$
such that  
$$
\lim_{i\to\infty}{\frac{1}{n_i\sigma(n_i)}}\log{a_{n_i}} =\infty
\ ,
\slabel{limsup}
$$
where $\sigma(n)$ is as in Proposition {5.2}.
Now we generate a sequence $g_0,g_1,\ldots,g_i,\ldots$ recursively, starting
with $g_0=f$, and taking, for all $i\ge 0$, $g_{i+1}=F(n_{i+1},g_i)$, where
$F(\cdot,\cdot)$ is as given in Proposition {5.2}. Each $g_i$ is a
$C^\infty$ critical circle map with $\rho(g_i)=\rho(f)$, and
$g_{i+1}=\Phi_{i+1}\circ g_i$, where $\Phi_{i+1}$ is a $C^\infty$ diffeo with 
$$
\|\Phi_{k+1}^{\pm 1} -
\roman{id}_{S^1}\|_{C^k}\;\le\;B_k\theta^{n_i(\sigma(n_i)-k+1)}  
\ ,
\slabel{Phi}
$$
for all $k$, where $0<\theta<1$ is a constant depending only on the real
bounds. From {\eqnoPhi} it follows that
$\Phi=\lim{\Phi_i\circ\cdots\circ\Phi_1}$ exists as a $C^\infty$ diffeo, and
therefore so does $g=\lim{g_i}=\Phi\circ f$ as a critical circle map.

Using properties ($c$) and ($d$) of Proposition {5.2} for each $g_i$, we
deduce that the $n_i$-th order discrepancy between $f$ and $g$ satisfies
$$
\left|\frac{|f_{n_i}^{N_i-1}(L_{n_i})|}{|f_{n_i}^{N_i-1}(R_{n_i})|} -
\frac{|g_{n_i}^{N_i-1}(\tilde{L}_{n_i})|}{|g_{n_i}^{N_i-1}(\tilde{R}_{n_i})|}
\right| \;\ge\; C|J_{n_i}|^{2\sigma(n_i)}
\ ,
\slabel{moredis}
$$
where $N_i=\lceil a_{n_i}/2\rceil$, etc.

Now, let $h:S^1\to S^1$ be the conjugacy between $f$ and $g$ mapping the
critical point $c$ to itself. Suppose $h$ were $C^{1+\beta}$ for some
$\beta>0$. Then the left-hand side of {\eqnomoredis} would be $\le
C|f_{n_i}^{N_i-1}(\Delta_1^{(n_i)})|^\beta$, where
$\Delta_1^{(n_i)}=L_{n_i}\cup R_{n_i}$. But by Yoccoz's Lemma, we have 
$$
|f_{n_i}^{N_i-1}(\Delta_1^{(n_i)})|\;\asymp\;\frac{1}{N_i^2}|J_{n_i}|\;\asymp\;
\frac{1}{a_{n_i}^2}|J_{n_i}|
\ .
\slabel{fjj}
$$
Combining the above with {\eqnomoredis} and {\eqnofjj}, we would get the
inequality 
$$
a_{n_i}^{2\beta}|J_{n_i}|^{2\sigma(n_i)-\beta}\;\le\;C
\ .
$$
But by the real bounds $|J_n|\ge C\mu^n$ for all $n$, where
$0<\mu<1$. Therefore, taking logarithms, we would have
$$
\frac{\beta\log{a_{n_i}}}{n_i\sigma(n_i)}\;\le\;\log{\frac{1}{\mu}} 
\ ,
$$
but this clearly contradicts {\eqnolimsup}. \qed

\heading Appendix A.\enspace
Compactness of renormalizations
\endheading

The real a-priori bounds proved in the section 3 have produced
a very important corollary, namely, that the renormalizations of 
an arbitrary $C^3$ critical circle map are uniformly bounded in the
$C^1$ topology. 
In this appendix we will use further a-priori estimates, this time involving
the Schwarzian derivative, to prove that such renormalizations are
uniformly bounded in the $C^{r-1}$ topology when the critical circle map is
$C^r$. Some technical tools are necessary.

\subhead 
A.1 The $C^m$-Approximation Lemma
\endsubhead 
In what follows, $m\ge 1$
will be a fixed integer and $I,J\subseteq {\Bbb R}$ fixed closed
intervals. We denote by $C^m(I)$ the Banach space of $C^m$-mappings
$f:I\to {\Bbb R}$ with the norm $\| f\|_{m}=\max\{\|D^if\|_0:\, 0\le
i\le m\}$, where $\|\phi\|_0=\sup_{x\in I}|\phi(x)|$. Sometimes, when
we need to emphasize the domain of $f$, we write $\|f\|_{I,m}$
instead of $\|f\|_m$.
We consider also the
closed, convex subset $C^m(I,J)\subseteq C^m(I)$ consisting of those
$f$'s such that $f(I)\subseteq J$.

Recall Leibnitz's formula for the $k$-th derivative of a product of
two functions,
$$
D^k(uv)\;=\;\sum_{j=0}^k {k\choose j}   D^ju\,D^{k-j}v\ \ ,
$$
from which it is clear that 
$$
\|uv\|_m\;\le\;2^m\|u\|_m\|v\|_m
\slabel{leib}
$$
whenever $u,v \in C^m(I)$. Something similar holds for the composition
of two $C^m$ mappings. Namely, we have Faa-di-Bruno's formula (cf.
{\refHea}, p.~42)
$$
D^k(f\circ g)\;=\;\sum_{j=1}^k B_{j,k}(D^1g,D^2g,\ldots ,D^jg)\,
D^{k-j+1}f\circ g \ , 
$$
where each $B_{j,k}$ is a homogeneous polynomial of degree $j$ on $j$
variables whose coefficients are non-negative numbers depending only
on $k$ and $j$. It readily follows from this formula that if $\psi\in
C^m(I,J)$ and $\phi\in C^m(J)$ then 
$$
\|\phi\circ\psi\|_m\;\le\;A(m)\|\phi\|_m\,\sum_{k=1}^m\|\psi\|_m^k \ ,
\slabel{faa}
$$
where $A(m)=\max_{1\le k\le m}\max_{1\le j\le k}\, B_{j,k}(1,1,\ldots ,1)$. 

Another well-known fact we will need below is the following ({\it
cf.\/} {\refFr}, Th.~3.1). Suppose
$m>1$ and consider the composition operator $(f,g)\mapsto f\circ g$ as
a map $\Theta :C^{m}(J)\times C^{m-1}(I,J)\to C^{m-1}(I)$ . Then
$\Theta$ is $C^1$ and its Fr\'echet derivative is given by  
$$
D\Theta(f,g)\,(u,v)\;=\;u\circ g+v\,Df\circ g \ .
\slabel{dcomp}
$$
Note that $C^m(J)\times C^{m-1}(I,J)\subseteq C^m(J)\times
C^{m-1}(I)$; we consider this last product endowed with the norm
$$
\big|(f,g)\big|_{I,J,m}=\max \{\|f\|_{J,m}, \|g\|_{I,m-1}\}
\ .
$$
\proclaim{Lemma {A.1}} For each $M>0$, there exists $c(M)>0$ such that, if
$f_1,g_1\in C^{m}(J)$ and $f_2,g_2\in C^{m-1}(I,J)$ and if
$|(f_1,f_2)|_{I,J,m} <M$ and $|(g_1,g_2)|_{I,J,m} <M$, then 
$$
\|f_1\circ f_2-g_1\circ g_2\|_{m-1}\;\le\;
c(M)\big|(f_1-g_1,f_2-g_2)\big|_{I,J,m} 
\ .
$$
\endproclaim
\demo{Proof} By the mean value theorem,
$$
\|f_1\circ f_2-g_1\circ g_2\|_{m-1}\;\le\; \sup_{(\phi
,\psi)}\|D\Theta(\phi ,\psi)\|\;|(f_1-g_1, f_2-g_2)|_{I,J,m} \ ,
$$
where the supremum is taken over all $(\phi ,\psi)$ in the line
segment joining $(f_1,f_2)$ to $(g_1,g_2)$ inside $C^{m}(J)\times
C^{m-1}(I,J)$, and where 
$$
\|D\Theta(\phi ,\psi)\|\;=\;\sup\big\{\|D\Theta(\phi
,\psi)(u,v)\|_{m-1}:\,|(u,v)|_{I,J,m}\le 1\big\}
$$
is the operator-norm of $D\Theta(\phi, \psi)$.
Using {\eqnodcomp}, and then {\eqnoleib} and {\eqnofaa}, we have
$$
\aligned
\|D\Theta(\phi ,\psi)(u,v)&\|_{m-1}\;\le \;\|u\circ\psi\|_{m-1} +
\|v\,D\phi\circ\psi\|_{m-1} \\
&\;\le \; \left(\|u\|_{m-1} +
2^{m-1}\|v\|_{m-1}\|D\phi\|_{m-1}\right)\,A(m-1)\,\sum_{k=1}^{m-1}
\|\psi\|_{m-1}^k 
\ .
\endaligned
$$
From this, and taking into account that $\|u\|_{m-1}\le \|u\|_{m}\le
|(u,v)|_{I,J,m}$ as well as $\|v\|_{m-1}\le |(u,v)|_{I,J,m}$, we
deduce that 
$$
\|D\Theta(\phi ,\psi )\|\;\le \;A(m-1)\,\left(1+2^{m-1}\|D\phi\|_{m-1}
\right)\, \sum_{k=1}^{m-1}\|\psi\|_{m-1}^k 
\ .
$$
Finally, since $\|D\phi\|_{m-1}\le\|\phi\|_{m}$ and $|(\phi ,\psi
)|_{I,J,m} < M$, we get
$$
\sup_{(\phi ,\psi )}\|D\Theta(\phi ,\psi )\|\;\le
\;A(m-1)\,\left(1+2^{m-1}M\right)\, \sum_{k=1}^{m-1} M^k \;=\;c(M)
\ .\qed
$$
\enddemo

Let us denote by ${\Bbb B}^m(I;M)$ the ball of radius $M$ centered at
the origin in $C^m(I)$.
\proclaim{Lemma {A.2}} {\sl{(The $C^m$-Approximation Lemma)}}

\noindent For each $M>0$, there exist constants $\varepsilon_M>0$ and 
$C_M>0$ such that the following holds for all $\varepsilon\le
\varepsilon_M$. Let $\Delta_1,\Delta_2,\ldots,\Delta_{n+1}$ be closed intervals on
the line or on the circle, and for each $1\le i\le n$ let $f_i,g_i\in
C^m(\Delta_i,\Delta_{i+1})$ be such that 
\itemitem{$(a)$} For all $1\le j\le k\le n$, we have $f_k\circ
f_{k-1}\circ\cdots\circ f_j\in {\Bbb B}^{m}(\Delta_j;M)$;
\itemitem{$(b)$} We have $\sum_{i=1}^n\|f_i-g_i\|_m<\varepsilon$.

\noindent Then for all $k\le n$ we have $g_k\circ
g_{k-1}\circ\cdots\circ g_1\in {\Bbb B}^{m-1}(\Delta_1;2M)$, and moreover 
$$ 
\|f_k\circ f_{k-1}\circ\cdots\circ f_1 -g_k\circ
g_{k-1}\circ\cdots\circ
g_1\|_{m-1}\;\le\;C_M\,\sum_{j=1}^k\|f_j-g_j\|_m
\ .
$$
\endproclaim
\demo{Proof} In the notation of Lemma {A.1}, let us write
$$
C_M\;=\;\max\left\{1,\, c(2M), \,c(2M)c(3M)\right\}
$$
and $\varepsilon_M=M/C_M$. We proceed by induction on $k$. When
$k=1$, we have $\|f_1-g_1\|_m\allowbreak\le \varepsilon$ and there is
nothing to prove. Suppose the assertion is valid for all $j<k$, and
write 
$$
\eqalign{\|f_k f_{k-1}\cdots f_1 - &g_k
g_{k-1}\cdots g_1\|_{m-1}\;\le\cr
&\sum_{j=1}^{k}\|f_k
\cdots  f_{j+1} g_j g_{j-1}\cdots g_1
-f_k\cdots f_{j+1} f_j g_{j-1}\cdots 
g_1\|_{m-1} \ .}
\slabel{fg}
$$
Since $|(f_j,\,g_{j-1}\circ\cdots \circ g_1)|_{\Delta_1,\Delta_j,m}< 2M$ and
also $|(g_j,\,g_{j-1}\circ\cdots \circ g_1)|_{\Delta_1,\Delta_j,m}< 2M$, it
follows from Lemma {A.1} that
$$
\|f_j  g_{j-1} \cdots   g_1- g_j   g_{j-1} \cdots
g_1\|_{m-1}\;\le \;c(2M)\|f_j-g_j\|_{m}
\ ,
$$
for $j=1,\ldots ,k$. In particular, by the induction hypothesis, we
have for all $1\le j\le k-1$
$$
\|f_j  g_{j-1} \cdots   g_1\|_{m-1}\;\le \;\|g_j 
g_{j-1} \cdots   g_1\|_{m-1}+\varepsilon_M \,c(2M)<3M
\ .
$$
Taking this back to {\eqnofg} and applying Lemma {A.1} again, we get
$$
\eqalign{
\|f_k  f_{k-1} \cdots  f_1 &-g_k 
g_{k-1} \cdots  g_1\|_{m-1}\cr
&\le\;c(2M)\|f_k-g_k\|_{m}+c(2M)c(3M)\,\sum_{j=1}^{k-1}\|f_j-g_j\|_{m}\;
\cr
&\le\;C_M\,\sum_{j=1}^k\|f_j-g_j\|_{m}\ \ ,}
$$
and this shows also that $\|g_k  g_{k-1} \cdots 
g_1\|_{m-1}\le M+\varepsilon_M C_M<2M$, thereby completing the induction.
\qed 
\enddemo

\subhead
A.2 Koebe principle revisited
\endsubhead
We present a generalization of the
classical Koebe non-linearity principle. This principle
states that if a $C^3$ diffeomorphism has non-negative Schwarzian
derivative on an open interval, then its non-linearity on any smaller
closed subinterval with {\it space} on both sides is bounded. The
generalized version below seems to be new.
We denote by $S\phi$ the Schwarzian derivative of $\phi$.

\proclaim{Lemma {A.3}} Given positive constants $B$ and $\tau$, there
exists $K_{\tau,B}>0$ such that the following holds. If $\phi$ is a
$C^3$-diffeomorphism of an interval $I\supseteq [-\tau,
1+\tau]$ into the reals and if $S\phi(t)\ge -B$ for all $t\in I$, 
then for all $t\in [0,1]$ we have 
$$
\left|\frac{\phi''(t)}{\phi'(t)}\right|\;\le\; K_{\tau,B}\ .
$$
\endproclaim

\demo{Proof} Writing $y=\phi''/\phi'$, so that $S\phi=y'-{1\over 2}y^2$, we
have the differential inequality
$$
y'\ge {1\over 2} y^2 - B\ .\slabel{ineqiff}
$$
Let $0\le t_0\le 1$ be a point where $|y(t)|$ attains its maximum in
$[0,1]$ and suppose $y_0=y(t_0)$ is such that $|y_0|>\sqrt{2B}=\beta$.
If $z(t)$ is the solution of the differential equation corresponding
to {\eqnoineqiff} with initial condition $z(t_0)=y_0$, then by a
well-known comparison theorem we must have $y(t)\ge z(t)$ for all
$t\ge t_0$ and $y(t)\le z(t)$ for all $t\le t_0$. Now, if $y_0>\beta$
then integration of the ODE leads to 
$$
z(t)=\beta \frac{(y_0+\beta)+(y_0-\beta)e^{\beta
(t-t_0)}}{(y_0+\beta)-(y_0-\beta)e^{\beta (t-t_0)}}\ . 
$$
Since this solution explodes at time 
$$
t_1=t_0 + \frac{1}{\beta}\log{\left(\frac{y_0+\beta}{y_0-\beta}\right)}
\ ,
$$
so does $y(t)$. Hence $t_1\notin I$, i.e. $t_1-t_0>\tau$, which gives us
$$
\frac{\phi''(t_0)}{\phi'(t_0)}=y_0< \beta\,
{{e^{\beta\tau}+1}\over{e^{\beta\tau}-1}}\ . 
$$
If instead $y_0<-\beta$, then we get 
$$
z(t)=\beta {{(\beta +y_0)-(\beta -y_0)e^{\beta (t-t_0)}}\over
{(\beta +y_0)+(\beta -y_0)e^{\beta (t-t_0)}}}\ ,
$$
and arguing as before for $t\le t_0$ gives us
$$
\frac{\phi''(t_0)}{\phi'(t_0)}=y_0 >
-\beta\,{{e^{\beta\tau}+1}\over{e^{\beta\tau}-1}}
\ . 
$$
Therefore the lemma is proved if we take
$$
K_{\tau,B}=\beta {{e^{\beta\tau}+1}\over{e^{\beta\tau}-1}}\ .\ \qed
$$
\enddemo
\subsubhead Remark\endsubsubhead As $B\to 0$, $K_{\tau,B}\to 2/\tau$ and
we recover the classical Koebe principle.

\subhead
A.3 Bounding the $C^2$ norms
\endsubhead
As before, let $f\in {\roman{Crit}}^r(S^1)$, $r\ge 3$, be a critical
circle map with critical point $c$ of power-law $p>1$. Conjugating
$f$ by a suitable $C^r$-diffeomorphism, we may assume that there
exists a neighborhood ${\Cal U}\subseteq {\Bbb R}/{\Bbb Z}$ of $c$
such that  
$$
f(x)\;=\;(x-c)|x-c|^{p-1}+a
$$
for all $x\in {\Cal U}$, where $a$ is a constant. This will be our
standing hypothesis on $f$, and we will sometimes say that $f$ is
a {\it canonical} circle map. Note in this case that for all 
$x\in {\Cal U}\setminus \{c\}$, the Schwarzian derivative of $f$
equals 
$$
Sf(x)\;=\;-\frac{p^2-1}{2(x-c)^2}\ .
\slabel{schwcr}
$$

We are interested in
the maps $f^{q_{n-1}-1}\colon I_n^1\to I_n^{q_{n-1}}$ and
$f^{q_n-1}\colon I_{n-1}^1\to I_{n-1}^{q_n}$, for a fixed $n\ge 1$.
They extend as diffeomorphisms to maximal open intervals
$J_{n,1}^-\supseteq I_n^1$ and $J_{n,1}^+\supseteq I_{n-1}^1$
respectively. 
When linearly rescaled to unit size, these diffeomorphisms are
called the {\it coefficients} of the $n$-th renormalization of $f$.

Let us be more precise. 
Consider the $n$-th renormalization of $f$, namely the commuting pair
$\bk(f)_n:[\lambda_n,1]\to {\Bbb R}$ defined in section 3.
We write
$J_{n,i}^-=f^{i-1}(J_{n,1}^-)$ for each $1\le i\le q=q_{n-1}$ and
$J_{n,j}^+=f^{j-1}(J_{n,1}^+)$ for each $1\le j\le Q=q_n$. We also write
$J_{n,0}^-=f^{-1}(J_{n,1}^-)$ and $J_{n,0}^+=f^{-1}(J_{n,1}^+)$.
For each $0\le j\le Q$, let $\Lambda_j\colon {\Bbb R}\to {\Bbb R}/{\Bbb
Z}$ be the affine (orientation-preserving) covering map such that
$\Lambda_j([0,1])= I_{n-1}^j$. Let $\Delta_n^-$ be the component of
$\Lambda_1^{-1}(J_{n,1}^-)$ that contains the interval
$[\lambda_n,0]$, and let $\Delta_n^+$ be the component of
$\Lambda_1^{-1}(J_{n,1}^+)$ that contains the interval $[0,1]$. Then 
define 
$$
\cases
{\Cal F}_n^- =\Lambda_0^{-1}\, f^{q-1}\,\Lambda_1 : \Delta_n^- 
\to {\Bbb R}&{}\cr
{}&{}\cr
{\Cal F}_n^+ = \Lambda_0^{-1}\, f^{Q-1}\,\Lambda_1 : \Delta_n^+\to
{\Bbb R}&{}\cr \endcases 
\ .
$$
These are the $n$-th renormalization coefficients of $f$. 
Consider also the so-called {\it folding factors} of $\bk(f)_n$,
namely the maps
$$
\cases
\varphi_n^- =\Lambda_1^{-1}\, f\,\Lambda_0 : \Lambda_0^{-1}(J_{n,0}^-)
\to {\Bbb R}&{}\cr 
{}&{}\cr
\varphi_n^+ =\Lambda_1^{-1}\, f\,\Lambda_0 : \Lambda_0^{-1}(J_{n,0}^+)
\to {\Bbb R}&{}\cr 
\endcases 
\ .
$$
Each of these maps is a homeomorphism with a unique critical point at
zero. One verifies at once that the maps $\bk(F)_n^-={\Cal F}_n^-\circ
\varphi_n^-$ and $\bk(F)_n^+={\Cal F}_n^+\circ \varphi_n^+$ are $C^r$
extensions of $\bk(f)_n^-$ and $\bk(f)_n^+$, respectively.

It will be useful to express the coefficients ${\Cal F}_n^\pm$ as long
compositions of rescaled diffeomorphisms in the following way. We will
give the explicit decomposition for ${\Cal F}_n^+$. 
A similar decomposition can be worked out for ${\Cal F}_n^-$.
Let us denote by $\Delta_{n,j}^+$ the component of
$\Lambda_j^{-1}(J_{n,j}^+)$ containing the unit interval. Note in
particular that $\Delta_n^+=\Delta_{n,1}^+$.
For each $j$ in the range $0\le j\le Q-1$, let 
$$
f_j\;=\;\Lambda_{j+1}^{-1}\, f\,\Lambda_j : \Delta_{n,j}^+\to
\Delta_{n,j+1}^+
\ .
$$
We call such maps the {\it elementary factors} of ${\Cal F}_n^+$.
Each $f_j$ is a $C^r$ diffeomorphism such that $f_j([0,1])=[0,1]$ (see
Figure 5). We have of course $\varphi_n^+=f_0$, but more importantly  
$$
{\Cal F}_n^+\;=\;\left(\Lambda_0^{-1}\circ \Lambda_Q\right)\circ
\left(f_{Q-1}\circ\cdots\circ f_j\circ\cdots f_1\right)
\ .
\slabel{rescomp}
$$
We note also that for all $t\in \Delta_{n,j}^+$ 
$$
Sf_j(t)\;=\;Sf\left(\Lambda_j(t)\right)\,\left[D\Lambda_j(t)\right]^2
\;=\;Sf\left(\Lambda_j(t)\right)\,|I_{n-1}^j|^2
\ ,
\slabel{reschw}
$$
by the chain rule for the Schwarzian derivative.

\demo{Notation} 
Given $J=[a,b]\subseteq {\Bbb R}$ and $\tau >0$,
we denote by $J^\tau$ the interval $[c,d]\supseteq J$ such that
$(a-c)/(b-a)=(d-b)/(b-a)=\tau$. Note that $J$ has {\it space} equal to
$\tau$ inside $J^\tau$.
\enddemo 

$$
\psannotate{
\psboxto(10.5cm;0cm){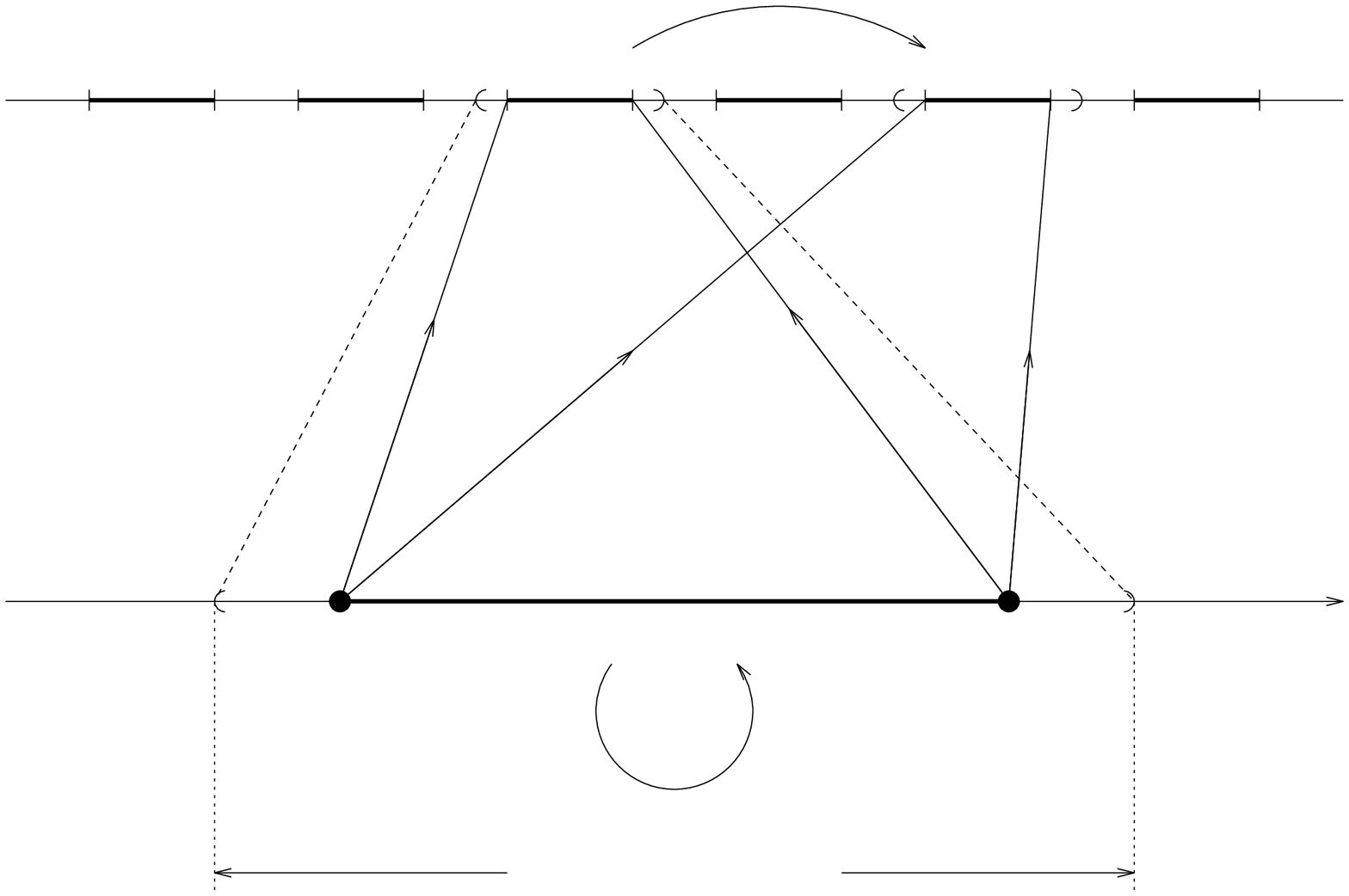}}
{\fillinggrid\at{4.1\pscm}{11.4\pscm}{$c$}
             \at{4.3\pscm}{12.5\pscm}{$I_{n-1}$}
             \at{10.5\pscm}{12.5\pscm}{$I_{n-1}^j$}
             \at{13.7\pscm}{14\pscm}{$f$}
             \at{16.4\pscm}{12.5\pscm}{$I_{n-1}^{j+1}$}
             \at{22\pscm}{11.4\pscm}{${\Bbb R}/{\Bbb Z}$}
             \at{8.1\pscm}{8\pscm}{$\Lambda_j$}
             \at{17.2\pscm}{8\pscm}{$\Lambda_{j+1}$}
             \at{6\pscm}{3.7\pscm}{$0$}
             \at{16\pscm}{3.7\pscm}{$1$}
             \at{20\pscm}{3.7\pscm}{${\Bbb R}$}
             \at{10.3\pscm}{2.5\pscm}{$f_j$}
             \at{10\pscm}{0.2\pscm}{$\Delta_{n,j}^+$}
}
$$
\bigskip
\centerline{\sl Figure 5. The elementary factors of ${\Cal F}_n^+$.}
\bigskip

\proclaim{Theorem {A.4}} ({\sl{The $C^2$ bounds}}) 

\noindent Let $f\in {\roman{Crit}}^3(S^1)$ be a
critical circle map with arbitrary irrational rotation number, let
$\bk(f)_n :[\lambda_n,1]\to {\Bbb R}$ be the $n$-th renormalization of
$f$, and let ${\Cal F}_n^{\pm}:\Delta_n^{\pm}\to {\Bbb R}$ be the
coefficients of $\bk(f)_n$. Also, let $f_j:\Delta_{n,j}^+\to
\Delta_{n,j+1}^+$ be the elementary factors of ${\Cal F}_n^+$.
There exist positive constants $B$ and
$\tau$ depending only on the real bounds for $f$ such that the
following statements hold for all $n\ge 1$.  
\itemitem{($a$)} We have $\Delta_n^-\supseteq [\lambda_n, 0]^\tau$ and 
$\Delta_n^+\supseteq [0, 1]^{\tau}=[-\tau, 1+\tau]$.
\itemitem{($b$)} For all $0\le j\le Q$, we have
$\Delta_{n,j}^+\supseteq [0, 1]^{\tau}$. 
\itemitem{($c$)} We have $|S{\Cal F}_n^-(t)|\le B$ for all $t\in
\Delta_n^-$ and $|S{\Cal F}_n^+(t)|\le B$ for all $t\in \Delta_n^+$. 
\itemitem{($d$)} More generally, for all $1\le j<k\le Q$, we have  
$|S(f_k\circ\cdots\circ f_j)(t)|\le B$ for all $t\in \Delta_{n,j}^+$. 
\itemitem{($e$)} The $C^2$ norms of the restrictions 
${\Cal F}_n^-|[\lambda_n,0]^{\tau/2}$ and ${\Cal
F}_n^+|[0,1]^{\tau/2}$ are bounded by $B$.
\itemitem{($f$)} More generally, for all $1\le j<k\le Q$, the $C^2$
norm of the restriction of $f_k\circ\cdots\circ f_j$ to
the interval $f_{j-1}\circ\cdots\circ f_1 ([0,1]^{\tau/2})$ is bounded by
$B$. 
\itemitem{($g$)} The $C^2$ norms of $\bk(f)_n^-$ and $\bk(f)_n^+$ are
bounded by $B$. 

\noindent Moreover, if $n$ is sufficiently large then both
coefficients have negative Schwarzian derivatives at all points of
their respective domains. 
\endproclaim

The proof will use the following lemma concerning the dynamical
partitions ${\Cal P}_n$.
Let us denote by $d(c,I)$ the distance between an interval
$I\subseteq S^1$ and the critical point $c$. For each $n\ge 1$, let  
$$
S_n\;=\;\sum_{I\in {\Cal P}_n\setminus
\{I_{n-1},I_n\}}\left(\frac{|I|}{d(c,I)}\right)^2 
\ .\slabel{Sn}
$$
\proclaim{Lemma {A.5}} The sequence $S_n$ is bounded (by a constant
depending only on $f$).
\endproclaim
\demo{Proof} Recall that ${\Cal P}_{n+2}$ is a strict refinement of
${\Cal P}_n$. From the real bounds, we know that there exists a
constant $0<\lambda<1$ depending only on $f$ such that, if $I$ is in
${\Cal P}_n$ and $J\subseteq I$ is in ${\Cal P}_{n+2}$, then $|J|\le
\lambda |I|$. Hence
$$
\sum_{I\supseteq J\in {\Cal P}_{n+2}}|J|^2\;\le
\;\left(\max_{I\supseteq J\in {\Cal P}_{n+2}}
|J|\right)|I|\;\le\;\lambda |I|^2  
\ .
$$
Since we also have $d(c,J)\ge d(c,I)$ whenever $J\subseteq I$, it
follows that
$$
\eqalign{
S_{n+2}\;\le\;\lambda S_n +& \sum_{{\Cal P}_{n+2}\ni J\subseteq
I_{n-1}\setminus I_{n+1}}\left(\frac{|J|}{d(c,J)}\right)^2 +  
\sum_{{\Cal P}_{n+2}\ni J\subseteq
I_n\setminus I_{n+2}}\left(\frac{|J|}{d(c,J)}\right)^2\cr
&{}\cr
&\le\;\lambda 
S_n+\lambda\left(\frac{|I_{n-1}|}{|I_{n+1}|}\right)^2+
\lambda\left(\frac{|I_n|}{|I_{n+2}|}\right)^2\ .}
$$
From this and the facts that $|I_{n-1}|\asymp |I_{n+1}|$ and
$|I_n|\asymp |I_{n+2}|$ , we get $S_{n+2}\le\allowbreak\lambda S_n+\mu$,
where $\mu$ is a constant depending only on $f$. But then, by
induction, 
$$
S_{2n}\;\le\;\lambda^{n-1} S_2+\frac{\mu}{1-\lambda}
\ \ ,\ \  
S_{2n+1}\;\le\;\lambda^n S_1+\frac{\mu}{1-\lambda}
\ ,
$$
and therefore $S_n$ is bounded as claimed. \qed
\enddemo

\demo{Proof of Theorem {A.4}} It is enough to prove this theorem
under the assumption that $f$ is canonical. The existence of $\tau>0$
such that ($a$) and ($b$) hold is a consequence of the real bounds. Hence
we proceed to prove ($c$) for ${\Cal F}_n^+$, the proof for ${\Cal
F}_n^-$ being completely similar. 
Making $\tau$ smaller if necessary and using the classical Koebe
non-linearity principle, we can assume that there exists $C>0$
depending only on the real bounds for $f$ such that 
$$
|D(f_j\cdots f_1)(t)|\le C
\ ,
\slabel{jone}
$$ 
for all $t\in [-\tau, 1+\tau]$ and all
$j=1,\ldots ,Q-1$. 

Let ${\Cal V}\subseteq S^1$ be an open set whose closure does not
contain $c$ and such that ${\Cal U}\cup {\Cal V}=S^1$. Also, let
$M=\sup_{x\in {\Cal V}}|Sf(x)|$.
We assume that
$n$ is so large that the largest interval in ${\Cal P}_n$ has length
smaller than the Lebesgue number of the covering $\{{\Cal U},{\Cal
V}\}$. 
Together with {\eqnorescomp} and {\eqnoreschw},
iterated use of the chain rule for the Schwarzian yields
$$
\eqalign{S{\Cal F}_n^+(t)\;&=\;S(f_{Q-1}\cdots f_j\cdots
f_1)(t)\cr
&=\;\sum_{j=1}^{Q-1}Sf_j(f_{j-1}\cdots
f_1(t))\,\left[D(f_{j-1}\cdots f_1)(t)\right]^2\cr
&=\;\sum_{j=1}^{Q-1}Sf(\Lambda_j\, f_{j-1}\cdots
f_1(t))\, |I_{n-1}^j|^2\,\left[D(f_{j-1}\cdots f_1)(t)\right]^2
\ .}
$$
We split this last sum into $\Sigma_1(t)+\Sigma_2(t)$, where
$\Sigma_1(t)$ is the sum over all $j$'s such that $I_{n-1}^j\subseteq
{\Cal U}$ and $\Sigma _2(t)$ is the sum over the remaining terms
(i.e. those with $I_{n-1}^j\subseteq {\Cal V}$). Then we have on one
hand 
$$
\left|\Sigma_2(t)\right|\;\le\;C^2M\sum_{I_{n-1}^j\subseteq {\Cal V}}
|I_{n-1}^j|^2\;\le\; C^2M\,\max_{1\le j\le Q-1}|I_{n-1}^j| 
\ .
\slabel{siga}
$$
On the other hand, since $d(c,J_{n,j}^+)\asymp d(c,I_{n-1}^j)$ for all
$j$, we have by {\eqnoschwcr} 
$$
\left|\Sigma_1(t)\right|\;\le\;C^2\sum_{I_{n-1}^j\subseteq {\Cal U}}
\frac{|I_{n-1}^j|^2}{[d(c,J_{n,j}^+)]^2}
\le C'S_n
\ ,\slabel{sigb}
$$
where $C'$ is another constant depending only on $f$ and $S_n$ is
given by {\eqnoSn}.
From {\eqnosiga} and {\eqnosigb} it follows that $|S{\Cal F}_n^+(t)|$ is
uniformly bounded, and this proves ($c$). Moreover, since by
{\eqnosiga} $\Sigma_2(t)$ goes to zero with $n$ while $\Sigma_1(t)$
is always negative and bounded away from zero, we deduce that
$S{\Cal F}_n^+(t)<0$ for all $n$ sufficiently large. The proof of
($d$) is entirely analogous.

To prove ($e$), let $B_0$ be the upper-bound that we have just
obtained for $|S{\Cal F}_n^+|$. Applying Lemma {A.3} to ${\Cal
F}_n^+$, we get for all $t\in [0,1]^{\tau/2}$
$$
\left|\frac{D^2{\Cal F}_n^+(t)}{D{\Cal F}_n^+(t)}\right|\;\le\;
K_{\tau_0,B_0}
\ ,
$$
where $\tau_0=\tau/2(1+\tau)$ is the space of $[0,1]^{\tau/2}$ inside
$[0,1]^\tau$. Therefore $\|D^2{\Cal F}_n^+\|_0\le\allowbreak
K_{\tau_0,B_0}\|D{\Cal F}_n^+\|_0\le CK_{\tau_0,B_0}$, by {\eqnojone}
above. 
This shows that the $C^2$ norm of ${\Cal F}_n^+$ is bounded as
claimed. 
A similar argument proves ($f$).
Finally, ($g$) follows from ($e$) and the fact that the
folding factors $\varphi_n^{\pm}$ are linear blow-ups of a fixed
power-law map. The theorem is therefore proved if we take $B$ to be
the largest of all the upper-bounds obtained in the argument.\qed
\enddemo

\demo{Remark}
We can go a bit further in ($e$), ($f$) and ($g$) and bound also the
$C^3$ norms. For this purpose, it suffices to note for instance that 
$$
D^3{\Cal F}_n^+(t)\;=\;D{\Cal F}_n^+(t)\,\left(S{\Cal F}_n^+(t)
+\frac{3}{2}\left[\frac{D^2{\Cal F}_n(t)}{D{\Cal F}_n(t)}\right]^2\right)
\ ,
$$
and then use ($c$) and ($e$). However, this argument does not
generalize to get bounds for higher derivatives. Our
bootstrap argument in the next section will follow a different route,
based on the $C^m$ Approximation Lemma.
\enddemo

\subhead 
A.4 Bounding the $C^{r-1}$ norms
\endsubhead
We will show that the sequence of renormalizations of a $C^r$ critical circle
map is bounded in the $C^{r-1}$ sense.
The limits fall into (a compact subset of) a special family
of analytic  critical commuting pairs known as the Epstein class. Moreover, we
will prove that such limits are attained at an exponential rate in the 
$C^{r-1}$ topology. The rate of convergence turns out to depend only
on the rotation number of the given critical circle map.  

An {\it Epstein map} is a homeomorphism $\varphi:I\to J$ between
closed intervals on the real line such that $\varphi^{-1}$ is the
restriction of an analytic univalent map $\Phi:{\Bbb C}(J')\to {\Bbb
C}(I')$, where $I'\supseteq I$ and $J'\supseteq J$ are open intervals.
Here we use the notation ${\Bbb C}(\Delta)=({\Bbb C}\setminus {\Bbb R})\cup
\Delta$. For example, every fractional linear transformation in
$PSL_2({\Bbb R})$ is an Epstein map when restricted to an interval on
the line which does not contain any of its poles. Further examples
include polynomial or rational diffeomorphisms with real coefficients.

\demo{Definition} A commuting pair $\bk(f)$ is said to be an {\it
Epstein} commuting pair if $\bk(f)^+=\allowbreak\varphi^+\circ Q$ and
$\bk(f)^-=\varphi^-\circ Q$, where $\varphi^+,\varphi^-$ are Epstein
maps and $Q$ is the power-law map $x\mapsto x^p$ (for some $p>1$). 
\enddemo 

\proclaim{Theorem {A.6}} 
Let $r\ge 3$ and let $f$ be a $C^r$
critical circle map with arbitrary irrational rotation number. Then
the sequence of renormalizations $\{{\Cal R}^n(f)\}$ is bounded in the
$C^{r-1}$ metric and converges $C^{r-1}$ exponentially fast to the
Epstein class. 
\endproclaim

The idea behind the proof
of Theorem {A.6} is quite simple. In the long composition defining
the $n$-th renormalization of a critical circle map, we replace the
factors away from the critical point by suitable fractional linear
approximations, which are all Epstein maps. The factors which are
close to the critical point are already Epstein because the map is
assumed to be a power-law there. Therefore the entire new composition
is an Epstein map. The Moebius approximations have to be carefully
chosen, however, so that the total error involved, estimated with the
help of the $C^m$ Approximation Lemma, be exponentially small in $n$
(the step of renormalization). We now present the technical result
which is needed. 

\proclaim{Lemma {A.7}} Given $r\ge 3$ and an orientation preserving
$C^r$-diffeomorphism $\phi:I\to {\Bbb R}$ of a closed interval $I$
onto its image, there exist constants $\ell_\phi>0$ and $K_\phi>0$
with the following property. For each 
closed interval $\Delta\subseteq I$ of length $|\Delta|\le \ell_\phi$, there
exists a fractional linear transformation $T_\Delta\in PSL_2({\Bbb R})$
with $T_\Delta(\Delta)=\phi(\Delta)$ such that, 
\itemitem{($a$)} $\sup_{x\in
\Delta}\left|D^k\phi(x)-D^kT_\Delta(x)\right|\;\le\;K_\phi|\Delta|^{3-k}$ for
$k=0,1,2$.
\itemitem{($b$)} $\sup_{x\in \Delta}\left|D^kT_\Delta(x)\right|\;\le\;K_\phi$
for all $1\le k\le r$.
\endproclaim
\demo{Proof} Let $\ell_\phi$ be the constant
$$
\ell_\phi\;=\;\min\left\{1,
\inf_{x\in I}\left|\frac{\phi'(x)}{\phi''(x)}\right|\right\} 
\ .
$$
Take any closed interval $\Delta\subseteq I$ with $|\Delta|\le
\ell_\phi$, and let $x_0$ be the left endpoint of $\Delta$. Let $T$ be
the unique fractional linear transformation with the same $2$-jet as
$\phi$ at $x_0$. Thus, if $T(x) = (a(x-x_0)+b)/(c(x-x_0)+d)$ with
$ad-bc=1$, then the coefficients are uniquely determined by the conditions
$$
\aligned
&T(x_0)\;=\;\frac{b}{d}\;=\;\phi(x_0)\ ,\\
&T'(x_0)\;=\;\frac{1}{d^2}\;=\;\phi'(x_0)\ , \\
&T''(x_0)\;=\;\frac{-2c}{d^3}\;=\;\phi''(x_0)\ .
\endaligned
\slabel{moba}
$$
Moreover, for all $k\ge 1$,
$$
D^kT(x)\;=\;\frac{(-1)^{k+1}k!c^{k-1}}{\left[c(x-x_0)+d\right]^{k+1}}
\ .
\slabel{mobb}
$$
Since $|x-x_0|\le |\Delta|\le\ell_\phi\le |\phi'(x_0)|/|\phi''(x_0)|=
|d|/2|c|$ for all $x\in \Delta$, we have 
$$
\frac{1}{2}|d|\;\le\;|c(x-x_0)+d|\;\le\;\frac{3}{2}|d|
\slabel{mobba}
$$ 
for each such $x$. Combining {\eqnomobb} with the lower bound in
{\eqnomobba}, we get
$$
\left|D^kT(x)\right|\;\le\;\frac{2^{k+1}k!|c|^{k-1}}{|d|^{k+1}}\;=\;
\frac{4k!|\phi''(x_0)|^{k-1}}{|\phi'(x_0)|^{k-2}}
\ ,
$$
for all $x\in \Delta$ and all $k\ge 1$, and consequently
$$
\sup_{x\in \Delta}\left|D^kT(x)\right|\;\le\;C_0\;=\;\max_{1\le k\le
r} \,\sup_{x\in I}
\left\{\frac{4k!|\phi''(x)|^{k-1}}{|\phi'(x)|^{k-2}}\right\}  
\ ,
\slabel{mobc}
$$
when $1\le k\le r$.
In particular, from
$$
D^2\phi(x)-D^2T(x)\;=\;\int_{x_0}^{x}D^3\phi(t)\,dt -
\int_{x_0}^{x}D^3T(t)\,dt 
\ ,
$$
we deduce that 
$$
\align
\left|D^2\phi(x)-D^2T(x)\right|\;&\le\;\|D^3\phi\|_0\,|x-x_0| + 
\frac{24|\phi''(x_0)|^2}{|\phi'(x_0)|}\,|x-x_0| \\
&\le\;\left(\|D^3\phi\|_0+
C_0\right)\,|\Delta|
\ ,
\endalign
$$
for all $x\in \Delta$. Integrating this inequality twice, using
{\eqnomoba}, we get
$$
\sup_{x\in \Delta}\left|D^k\phi(x)-D^kT(x)\right|\;\le\;C_1|\Delta|^{3-k}
\ ,
\slabel{mobd}
$$
for $k=0,1,2$, where $C_1=C_0+\|D^3\phi\|_0$. 

Looking at {\eqnomobc} and {\eqnomobd}, we see that $T$ is almost what
we want, but not quite because in general it does not map $\Delta$ onto
$\phi(\Delta)$. To correct this flaw, we replace $T$ by $T_\Delta=A\circ T$,
where $A$ is the unique affine, orientation-preserving map that
carries $T(\Delta)$ onto $\phi(\Delta)$. We have
$$
A(t)-t\;=\;\left[\frac{|\phi(\Delta)|}{|T(\Delta)|}-1\right]\,
\left(t-T(x_0)\right) 
\ ,
\slabel{mobe}
$$
for all $t\in T(\Delta)$, because $\phi(x_0)=T(x_0)$. Let
$\mu=|\phi(\Delta)|/|T(\Delta)|$. Since by {\eqnomobd} we have
$\left||\phi(\Delta)|-|T(\Delta)|\right|\le 2C_1|\Delta|^3$, and 
since by the upper-bound in {\eqnomobba} we have 
$$
\frac{|T(\Delta)|}{|\Delta|}\;\ge\;\inf_{x\in \Delta}{\frac{1}{[c(x-x_0)+d]^2}}
\;\ge\; \frac{4}{9d^2}\;=\;\frac{4}{9}\,|\phi'(x_0)|
\ ,
$$
it follows that
$$
|\mu-1|\;\le\;\frac{9C_1}{2|\phi'(x_0)|}\,|\Delta|^2 
\ .
$$
Thus we see that, for all $t\in T(I)$,
$$
|A'(t)-1|\;=\;|\mu-1|\;\le\;\frac{9C_1}{2\inf_{x\in I}|\phi'(x)|}
\,|\Delta|^2\;=\;C_2|\Delta|^2  
\ ,
$$
On the other hand, since $|T(\Delta)|\le
\|D\phi\|_0|\Delta|+2C_1|\Delta|^3$, and since $|\Delta|\le 1$, it
follows from {\eqnomobe} that
$$
|A(t)-t|\;\le\;C_2\,\left(\|D\phi\|_0+2C_1\right)\,
|\Delta|^3\;=\;C_3|\Delta|^3
\ .
$$
Therefore
$$
|\phi(x)-T_\Delta(x)|\;\le\;|\phi(x)-T(x)|+|T(x)-A(T(x))|\;\le\;
(C_1+C_3)\,|\Delta|^3
\ ,
\slabel{mobfa}
$$
and moreover, using the fact that $D^kT_\Delta(x)= \mu D^kT(x)$ for all
$k$, 
$$
\aligned
|D^k\phi(x)-D^kT_\Delta(x)|\;&\le\;|D^k\phi(x)-D^kT(x)|+
|\mu-1|\,|D^kT(x)| \\   
&\le\;C_1|\Delta|^{3-k}+C_0C_2|\Delta|^2 \\
&\le\;(C_1+C_0C_2)\,|\Delta|^{3-k}\ ,
\endaligned
\slabel{mobg}
$$
for all $x\in \Delta$ and $k=1,2$. Finally, for all $k\ge 1$ we have
$$
|D^kT_\Delta(x)|\;\le\;\left(1+C_2|\Delta|^2\right)\,|D^kT(x)|\;\le\;
(1+C_2)\,C_0 
\ .
\slabel{mobh}
$$
Part $(a)$  now follows from {\eqnomobfa} and {\eqnomobg}, while part
$(b)$ follows from {\eqnomobh}, provided we take $K_\phi=\max\{C_1+C_3,
C_1+C_0C_2,(1+C_2)C_0\}$.\qed 
\enddemo

\demo{Proof of Theorem {A.6}}
We now expand the outline given above and present a complete proof of
Theorem {A.6}. In the proof, we will denote by $C_0, C_1,\ldots$
positive constants depending only on the real bounds for $f$.
As before, we may assume from the start that $f$ is {\it canonical},
and accordingly we consider the covering $\{{\Cal U}, {\Cal V}\}$ of
$S^1$ defined in the proof of Theorem {A.4}.
Since the folding factors of $\bk(f)_n$ are power-law maps, and
therefore already Epstein, it suffices to prove that the {\it
coefficients} of $\bk(f)_n$ can be approximated by Epstein maps, up to
an error exponentially small in $n$ in the $C^{r-1}$ topology. We will
do this for ${\Cal F}_n^+$, the proof for ${\Cal F}_n^-$ being the
same. 

As in the previous section, let
$f_j : \Delta_{n,j}^+\to\Delta_{n,j+1}^+$, $1\le j\le Q-1$, be the
elementary factors of ${\Cal F}_n^+$. For each $1\le j\le Q$ we define
$$
\Delta_j\;=\;f_{j-1}\circ\cdots\circ f_2\circ
f_1\left([0,1]^{\tau/4}\right)\;\subseteq\;\Delta_{n,j}^+
\ ,
$$
where $\tau$ is the constant of Theorem {A.4}. Note that
$f_j(\Delta_j)= \Delta_{j+1}$. Let
$\Delta_j'=\Lambda_j(\Delta_j)$, and observe also that
$I_{n-1}^j\subseteq\Delta_j'\subseteq J_{n,j}^+$. 

We introduce individual Epstein approximations $g_j$ to each $f_j$.
There are two cases to consider. It may happen that
$\Delta_j'\subseteq {\Cal U}$, in which case we simply take $g_j=f_j$.
Otherwise, we have $\Delta_j'\subseteq {\Cal V}$. In this case, we let
$T_j:\Delta_j'\to \Delta_{j+1}'$ be the Moebius approximation to
$f|\Delta_j'$ that we get applying Lemma {A.7} to the restriction of $f$ 
to ${\Cal V}$, and then take $g_j=\Lambda_{j+1}^{-1}\circ T_j\circ
\Lambda_j$. Note that $g_j(\Delta_j)=\Delta_{j+1}$.
\smallskip
\noindent
{\smc Claim 1.} {\it We have $\|f_j-g_j\|_r\le C_0|I_{n-1}^j|^2$ for
all $j$.}
\smallskip
This is obvious when  $I_{n-1}^j\subseteq {\Cal U}$. When
$I_{n-1}^j\subseteq {\Cal V}$, we have $|I_{n-1}^j|\asymp
|I_{n-1}^{j+1}|$, because the derivative of $f$ on ${\Cal V}$ is
bounded away from zero, and we also have $|\Delta_j'|\asymp
|I_{n-1}^j|$. 
Moreover, for all $1\le s\le r$ and all $x\in \Delta_j$,   
$$
D^sf_j(x)-D^sg_j(x)\;=\;\frac{|I_{n-1}^j|^s}{|I_{n-1}^{j+1}|}\,
\left(D^sf(\Lambda_j(x))-D^sT_j(\Lambda_j(x))\right)
\ .
$$
Therefore the claim follows from Lemma {A.7} (treat the cases $s=1,2$ 
separately). 

Now, recall from Theorem {A.4} that for all $1\le j<k\le Q-1$ we have
$$
\|f_k\circ\cdots\circ f_j\|_2\le B
\ .
$$
\smallskip
\noindent
{\smc Claim 2.} {\it If $n$ is sufficiently large then for all $1\le
j<k\le Q-1$ we have}
$$
\|f_k\circ\cdots\circ f_j-g_k\circ\cdots\circ g_j\|_1\;\le\;
C_1\max_{0\le i\le Q}|I_{n-1}^i|
\ .
\slabel{jkbounds}
$$

Take $n_0$ so large that $C_0\max{|I_{n_0-1}^j|}<\varepsilon_B$,
where $\varepsilon_B$ is the constant given by Lemma {A.2} when we take
$M=B$. Then from Claim 1 and {\eqnojkbounds}, the hypotheses of Lemma
{A.2} are satisfied, and we get for all $n\ge n_0$
$$
\eqalign{
\|f_k\circ\cdots\circ f_j-g_k\circ\cdots\circ g_j\|_1\;&\le\;
C_B\sum_{i=j}^k\|f_i-g_i\|_2\;\le\;C_0C_B\sum_{i=j}^k|I_{n-1}^i|^2\cr
&{}\cr
&\le\;C_0C_B\max_{0\le i\le Q}{|I_{n-1}^i|}
\ ,}
$$
where $C_B$ is the constant of Lemma {A.2} for $M=B$. This proves
the claim. 

In order to bootstrap these $C^1$ estimates up to $C^{r-1}$ estimates,
we apply the $C^m$ Approximation Lemma once more, this time reversing
the roles of $f_j$ and $g_j$, and with $m=r$. Thus, we need to verify the
hypotheses of that lemma in this new situation.
\smallskip
\noindent
{\smc Claim 3.} {\it For all $1\le j<k\le Q-1$, we have
$\|g_k\circ\cdots\circ g_j\|_r\le C_2$.}
\smallskip
For brevity, write $G_{jk}=g_k\circ\cdots\circ g_j$. Then
$G_{jk}^{-1}$ is univalent on ${\Bbb C}(\Delta_{jk})$, where
$\Delta_{jk}$ is an interval containing $G_{jk}(\Delta_j)$ with
definite space on both sides, by our choice of $\tau$. Using Koebe's
one-quarter theorem, it is not difficult to see that the domain
$\Omega_{jk}=G_{jk}^{-1}({\Bbb C}(\Delta_{jk}))$ contains a rectangle
$W_j=\Delta_j^{\alpha}\times [-\beta,\beta]$, and that $d(\partial
W_j, \partial\Omega_{jk})\ge \gamma$ where $\alpha$, $\beta$ and
$\gamma$ are positive constants depending only on $\tau$ and the real
bounds for $f$. Hence, from the complex Koebe's distortion theorem,
we get  
$$
\Big{|}\frac{G_{jk}'(z)}{G_{jk}'(w)}\Big{|}\;\le\;
\exp\left\{\frac{2}{\gamma}\,{\roman{diam}}\,(W_j)\right\}
\;\le\;C_3
\ ,
$$
for all $z,w \in W_j$. This together with the mean-value theorem
gives us $|G_{jk}'(z)|\le C_4$, and therefore also $|G_{jk}(z)|\le
C_5$, for all $x\in W_j$. Now we use Cauchy's integral formula to
bound all higher derivatives of $G_{jk}$. We have for all $x\in
\Delta_j$ and all $s\ge 1$
$$
|D^sG_{jk}(x)|\;=\;\frac{s!}{2\pi}\, \Big{|}\int_{\partial
W_j}\frac{G_{jk}(z)}{(z-x)^{s+1}}\,dz
\Big{|}\;\le\;\frac{C_5s!}{\pi}(\beta+(1+2\alpha)|\Delta_j|)\,
\delta_j^{-s-1}  
\ ,
$$
where $\delta_j=\inf_{x\in\Delta_j}{d(x,\partial
W_j)}=\min\{\alpha|\Delta_j|,\beta\}\ge \delta=\min\{a,b\}$.
Therefore $|D^sG_{jk}(x)|\le C_6s!\delta^{-s-1}$. This shows that
$\|G_{jk}\|_r$ is bounded as claimed.

From Claims 1 and 3, the hypotheses of Lemma {A.2} are therefore
satisfied, and we have 
$$
\eqalign{
\|f_k\circ\cdots\circ f_j-g_k\circ\cdots\circ g_j\|_{r-1}\;&\le\;
C_{C_2}\sum_{i=j}^k\|f_i-g_i\|_r\;\le\;C_0C_{C_2}\sum_{i=j}^k|I_{n-1}^i|^2\cr
&{}\cr
&\le\;C_0C_{C_2}\max_{0\le i\le Q}{|I_{n-1}^i|}
\ ,}
$$
this time for all $n$ large enough so that
$C_0\max{|I_{n-1}^j|}<\varepsilon_{C_2}$, where $C_{C_2}$ and
$\varepsilon_{C_2}$ are the constants of Lemma {A.2} for $M=C_2$. Since
$\max{|I_{n-1}^j|}$ decreases exponentially with $n$, we are done. \qed
\enddemo

\head Appendix B.\enspace
Proof of Yoccoz's Lemma\endhead

The main geometric idea behind the proof of Yoccoz's Lemma is to use
the negative Schwarzian property of $f$ to {\it squeeze} the graph of
$f$ between the graphs of two Moebius transformations. The required
estimate for $f$ will then follow from the corresponding estimate for
Moebius transformations, which we now state and prove.

Consider the fractional linear transformation $T(x)=x/(1+x)$, and
given $\varepsilon>0$, let $T_\varepsilon(x)=T(x)-\varepsilon$. We are
interested in certain quantitative aspects of the orbit
$x_n=T^n_\varepsilon (x_0)$ for $x_0=1$. Observe that this sequence is
strictly decreasing.
\proclaim{Lemma {B.1}} Let $N>0$ be such that $x_{N+1}\le 0<x_N$.
Then we have $N\asymp 1/\sqrt{\varepsilon}$ and moreover
$x_n-x_{n+1}\asymp 1/n^2$ for $n=0,1,\ldots ,N$.
\endproclaim
\demo{Proof} Writing $\delta_n=T^n(x_0)-T^n_\varepsilon (x_0)$,
we have
$$
\delta_n=\varepsilon +
\frac{\delta_{n-1}}{(1+\frac{1}{n})(1+\frac{1}{n}-\delta_{n-1})} 
\slabel{delta}
$$
for all $n=1,2,\ldots ,N+1$. We claim that
$$
\frac{n\varepsilon}{6}\;\le\; \delta_n\;\le\;n\varepsilon
\ .
\slabel{neps}
$$
The last inequality is clear. To prove the first, we note from
{\eqnodelta} that
$$
\delta_n\;\ge\;\varepsilon + \(\frac{n}{n+1}\)^2\delta_{n-1}
\ .
$$
By induction, this gives us
$$
\delta_n\;\ge\;\frac{\varepsilon}{(n+1)^2}\,\(1^2+2^2+\cdots +n^2\)
\;=\;\frac{\varepsilon}{(n+1)^2}\,\frac{n(n+1)(2n+1)}{6}\;\ge\;
\frac{n\varepsilon}{6} 
\ ,
$$
which proves the claim. Now, from the fact that $x_{N+1}\le 0<x_N$ we
have the inequalities 
$$
\delta_N\;<\;\frac{1}{N+1}\ , \ \ \delta_{N+1}\;\ge\;\frac{1}{N+2}\ 
\ .
$$
Then, using {\eqnoneps}, we get
$$
\frac{1}{(N+1)(N+2)}\;\le\;\varepsilon\;<\;\frac{6}{N(N+1)}
\ ,
\slabel{epsize}
$$
which proves the first assertion. 

Next, note that since
$\[x_{N+1},x_N\]\subseteq
\[T_\varepsilon(0),T_\varepsilon^{-1}(0)\]= \[-\varepsilon,
\varepsilon/(1-\varepsilon)\]$, we have
$$
\varepsilon\;<\;x_N-x_{N+1}\;<\;3\varepsilon
\slabel{delN}
$$
Hence, by {\eqnoepsize}, we get $x_N-x_{N+1}\asymp 1/N^2$ and the
second assertion is proved when $n=N$. To prove it in general using
this information, observe that
$$
x_n-x_{n+1}\;=\;\frac{x_{n-1}-x_n}{(1+x_{n-1})(1+x_n)}\;=\;
\frac{x_{n-1}-x_n}{(1+\frac{1}{n}-\delta_{n-1})(1+\frac{1}{n+1}-\delta_n)}
$$
implies
$$
x_n-x_{n+1}\;\ge\;\frac{n}{n+2}\,(x_{n-1}-x_n)
\ .
$$
By induction, this gives on one hand
$$
x_n-x_{n+1}\;\ge\;\frac{2}{(n+1)(n+2)}\,(x_0-x_1)\;\ge\;
\frac{1}{(n+1)(n+2)} 
\ ,
$$
and on the other hand, using {\eqnoepsize} and {\eqnodelN},
$$
x_n-x_{n+1}\;\le\;(x_N-x_{N+1})\,
\prod_{j=1}^{N-n}\left(\frac{n+j+2}{n+j}\right) 
\;<\;\frac{54}{(n+1)(n+2)} 
\ .
$$
This proves the second assertion in all cases. \qed
\enddemo

Now recall that $f:\Delta_1\cup\Delta_2\cup\cdots\cup\Delta_a\to {\Bbb
R}$ satisfies $f(\Delta_j)=\Delta_{j+1}$ for all $j$. Without loss of
generality, we can assume that $f(x)<x$ for all $x$. Thus, if we call
$x_0$ the right endpoint of $\Delta_1$ and write $x_j=f^j(x_0)$, we
have $\Delta_j=[x_j,x_{j-1}]$ for all $j$. Since our map $f$ is a
negative-Schwarzian diffeomorphism, there exists a unique
$z$ in the domain of $f$ such that $\varepsilon=|f(z)-z|\le |f(x)-x|$
for all $x$. Since the statement we want to prove is invariant under
affine changes of coordinates, we may assume also that $z=0$ and
$x_0=1$. In this setting, we want to prove that $|\Delta_j|\asymp
1/j^2$ for all $j$ such that $\Delta_j\subseteq [0,1]$. Note that
$f'(0)=1$. 

Next, let $A$ be the Moebius transformation on the line such that
$A(x_0)=f(x_0)$ and $A(0)=f(0)$ and $A'(0)=f'(0)=1$. This determines
$A$ uniquely, and in fact
$$
A(x)\;=\;\frac{x}{1+\lambda x} -\varepsilon
\ ,
$$
for some $\lambda>0$. Since $Sf<0$, we see that $A(x)\le f(x)$ for all
$x\in [0,1]$. 

Likewise, let $B$ be the Moebius transformation such that
$B(x_a)=f(x_a)$, $B(0)=f(0)$ and $B'(0)=f'(0)=1$. This determines $B$
uniquely, and in fact 
$$
B(x)\;=\;\frac{x}{1+\mu x} -\varepsilon
\ ,
$$
for some $\mu>0$. This time, since $x_a<0$ and $Sf<0$, we have
$f(x)\le B(x)$ for all $x\in [0,1]$. In particular, $\lambda>\mu$. It
is easy to see that $\lambda/\mu\le c_\sigma$, where $c_\sigma$
depends only on the constant $\sigma$ in the statement. 

\proclaim{Lemma {B.2}} Let $x\in [0,1]$ and $k>0$ be such that
$A(x)<B^k(x)$. Then $k\le 1+\lambda/\mu$.
\endproclaim
\demo{Proof} By induction we have 
$$
B^k(x)\;\le\;\frac{x}{1+(k-1)\mu x}-\varepsilon
\ .
$$
Therefore $A(x)<B^k(x)$ implies $(k-1)\mu x<\lambda x$. \qed
\enddemo

Now, let us write $\alpha_n=A^n(x_0)$ and $\beta_n=B(x_0)$. By Lemma
{B.2}, the number of $\beta_j$'s inside each interval of the form
$[\alpha_{n+1},\alpha_n]$ is bounded independently of $n$. Moreover,
since $\alpha_n<x_n<\beta_n$ for all $n$, the number of $x_j$'s inside
each $[\alpha_{n+1},\alpha_n]$ is also bounded independently of
$n$. To prove that $|\Delta_j|\asymp 1/j^2$, we proceed as 
follows. Let $\ell>0$ be such that $\beta_{\ell+1}\le x_j\le
\beta_{\ell}\le x_{j-1}$. Then Lemma {B.2} says that $\ell\le Cj$,
and we have also
$$
|\beta_{\ell+1}-\beta_{\ell}|\;<\;|B(x_{j-1})-x_{j-1}|\;<\;|x_j-x_{j-1}|
\ .
$$
Since by Lemma {B.1} we have 
$$
|\beta_{\ell+1}-\beta_{\ell}|\;\asymp\;\frac{1}{\ell^2}\;\ge\;\frac{1}{Cj^2}
\ ,
$$
it follows that $|\Delta_j|=|x_j-x_{j-1}|\ge 1/Cj^2$.

To prove an inequality in the opposite direction, let $m$ be the
largest integer such that $\alpha_m>x_{j-1}$. Then, again by Lemma
{B.2}, we have $j\le Cm$. Since $A(x)<f(x)<x$ for all $x$, we also
have $\Delta_j\subseteq [\alpha_{m+2},\alpha_m]$. Using Lemma
{B.1} once more, we deduce that
$$
|\Delta_j|\;\le\;\frac{C}{m^2}\;\le\;\frac{C}{j^2}
\ .
$$
This completes the proof of Yoccoz's Lemma. \qed



\Refs
\widestnumber\key{MMM8}

\ref \key\refnoBi
\by \Bi
\book Probability and measure
\publ John Wiley \& Sons
\yr 1986
\publaddr New York
\endref

\ref \key\refnoCa
\by \Ca
\paper On mappings conformal at the boundary
\yr 1967
\vol 19
\pages 1--13
\jour \JAM
\endref

\ref \key\refnodF
\by \dF
\paper Asymptotic rigidity of scaling ratios for critical circle
mappings 
\publ IMS Stony Brook preprint 96/13 (1996)
\jour to appear in \ETDS
\endref

\ref \key\refnodFdM
\by \dF\ and \Me
\paper Rigidity of critical circle mappings II
\publ IMS Stony Brook preprint 97/17 (1997)
\endref

\ref \key\refnoFr
\by \Fr
\paper Manifolds of $C^r$ mappings and applications to differentiable
dynamics 
\inbook Studies in Analysis, Advances in Mathematics Supplementary
Studies
\yr 1979
\vol 4
\pages 271--290
\publ Academic Press
\endref


\ref \key\refnoGd
\by \Gd
\paper The fractional dimension theory of continued fractions 
\yr 1941
\vol 37
\pages 199--228
\jour \PCPS
\endref


\ref \key\refnoGk
\by \Gk
\paper Harmonic scalings for smooth families of diffeomorphisms of the
circle 
\yr 1990
\vol 4
\pages 935--954
\jour \NL
\endref

\ref \key\refnoGkone
\bysame
\paper H\"older classes for circle maps
\jour to appear in \PAMS
\endref



\ref \key\refnoGS
\by \Gk\  and \Sw
\paper Critical circle maps near bifurcation
\yr 1996
\vol 176
\pages 227--260
\jour \CMP
\endref


\ref \key\refnoHea
\by \He
\paper Sur la conjugaison differentiable des diff\'eomorphismes du
cercle a des rotations 
\yr 1979
\vol 49
\pages 5--234
\jour \IHES
\endref

\ref \key\refnoHeb
\bysame
\paper Conjugaison quasi-sim\'etrique des hom\'eomorphismes du
cercle a des rotations 
\publ manuscript
\yr 1988
\endref

\ref \key\refnoLa
\by \La
\paper Renormalization group methods for critical circle mappings with
general rotation number
\inbook VIIIth International Congress on Mathematical Physics
(Marseille, 1986)
\publ World Sci. Publishing
\publaddr Singapore
\yr 1987
\pages 532--536
\endref

\ref \key\refnoLy
\by \Ly
\paper Combinatorics, geometry and attractors of quasi-quadratic mappings 
\yr 1994
\vol 140
\pages 347--404
\jour \AnM
\endref

\ref \key\refnoMc
\by \Mc
\jour Annals of Math Studies
\vol 142
\paper Renormalization and 3-manifolds which fiber over the circle
\publ Princeton University Press
\publaddr Princeton
\yr 1996
\endref

\ref \key\refnoMS
\by \Me \ and \vS
\book One dimensional dynamics
\publ Springer-Verlag
\yr 1993
\publaddr Berlin and New York
\endref

\ref \key\refnoRa
\by \Ra
\paper Global phase-space universality, smooth conjugacies and
renormalization: the $C^{1+\alpha}$ case 
\yr 1988
\vol 1
\pages 181--202
\jour \NL
\endref

\ref \key\refnoSwa
\by \Sw
\paper Rational rotation numbers for maps of the circle
\yr 1988
\vol 119
\pages 109--128
\jour \CMP
\endref

\ref \key\refnoSwb
\bysame
\paper One-dimensional maps and Poincar\'e metric
\yr 1992
\vol 4
\pages 81--108
\jour \NL
\endref

\ref \key\refnoYoa
\by \Yo
\paper Conjugaison differentiable des 
diffeomorphismes du cercle dont le nombre de rotation verifie une
condition Diophantienne
\yr 1984
\vol 17
\pages 333--361
\jour  Ann. Sci. de l'Ec. Norm. Sup. 
\endref

\ref \key\refnoYob
\bysame
\paper Il n'y a pas de contre-example de Denjoy 
analytique, 
\yr 1984
\vol 298
\pages 141--144
\jour C. R. Acad. Sci. Paris 
\endref

\endRefs
\enddocument